\documentclass{bmvc2k}
\usepackage{amsmath,amsthm,amssymb,bbm,MnSymbol,wasysym}
\graphicspath{{./Figures/}}

\title{Second order elastic metrics on the shape space of curves}

\addauthor{Martin Bauer}{bauer.martin@univie.ac.at}{1}
\addauthor{Martins Bruveris}{martins.bruveris@brunel.ac.uk}{2}
\addauthor{Philipp Harms}{philipp.harms@math.ethz.ch}{3}
\addauthor{Jakob M\o ller-Andersen}{jakmo@dtu.dk}{4}

\addinstitution{Technical University of Vienna}
\addinstitution{Brunel University London}
\addinstitution{ETH Z\"urich}
\addinstitution{Technical University of Denmark}

\runninghead{Bauer et al.}{Elastic metrics on the space of curves}

\def\th{\theta}

\def\ph{\varphi}

\let\on=\operatorname
\let\mc=\mathcal

\def\R{\mathbb{R}}

\theoremstyle{plain}
\newtheorem{theorem}[subsection]{Theorem}

\newtheorem*{theorem*}{Theorem}
\newtheorem*{proposition*}{Proposition}
\newtheorem*{lemma*}{Lemma}

\theoremstyle{definition}
\newtheorem{definition}[subsection]{Definition}
\newtheorem{remark}[subsection]{Remark}

\newtheorem*{definition*}{Definition}
\newtheorem*{remark*}{Remark}
\newtheorem*{example*}{Example}
\newtheorem*{openquestion*}{Open Question}

\begin{document}

\maketitle
\begin{abstract}
Second order Sobolev metrics on the space of regular unparametrized planar curves have several desirable completeness properties not present in lower order metrics, but numerics are still largely missing. In this paper, we present algorithms to numerically solve the initial and boundary value problems for geodesics. The combination of these algorithms allows to compute Karcher means in a Riemannian gradient-based optimization scheme. Our framework has the advantage that the constants determining the weights of the zero, first, and second order terms of the metric can be chosen freely. Moreover, due to its generality, it could be applied to more general spaces of mapping. We demonstrate the effectiveness of our approach by analyzing a collection of shapes representing physical objects.
\end{abstract}

\section{Introduction}
\label{sec:intro}
 
\footnotetext{\dagger All authors contributed equally to the article}

Unparametrized curves arise naturally in shape analysis and its many applications, including medical imaging \cite{Younes2012,Pennec2015gsi,XKS2014}, object tracking \cite{Sundaramoorthi2008,Sundaramoorthi2011}, computer animation \cite{Esl2014,Esl2014b}, computer aided design \cite{kmp07}, speech recognition \cite{su2014b}, analysis of bird migration patterns  and hurricane paths \cite{su2014}, biology \cite{Dryden1998,Laga2014}, and many other fields \cite{Krim2006,Bauer2014}. In many instances, the rationale for identifying curves differing only by a reparametrization is that the curves represent the boundaries of physical shapes. Shapes can be analyzed mathematically by endowing the space of shapes with a Riemannian metric. Riemannian shape analysis has developed into an  active field of research by now.

To do statistics on shape space, robust and efficient implementations of the boundary and initial value problems for geodesics are needed. Unfortunately, the arguably simplest metric on shape space is degenerate: it is well-known that the Riemannian distance induced by the $L^2$-metric vanishes on spaces of parametrized and unparametrized curves \cite{Michor2006c}. The discovery of this degeneracy led to an investigation of Sobolev metrics of higher order \cite{Michor2007,Younes1998,Shah2013}. 

For the important class of first order metrics numerics are well developed by now, but a restriction on the parameters of the metric must be imposed \cite{Klassen2004,Jermyn2011}. The situation is drastically different for second order metrics, despite the fact that they enjoy better completeness properties \cite{Bruveris2014,Bruveris2014b_preprint}. On the positive side, the geodesic boundary value problem under second order Finsler metrics on the space of $BV^2$ curves was implemented numerically in \cite{Vialard2014_preprint}. Moreover, on the space of parametrized curves, there are numerics for the boundary and initial value problems for geodesics under second order Sobolev metrics \cite{BBHM2015,Bauer2014c}. On spaces of unparametrized curves numerics for second order Sobolev metrics are, however, still lacking. This is an important challenge and the topic of this paper. 

We present a numerical implementation of the initial and boundary value problem for geodesics of planar, unparametrized curves under second order Sobolev metrics. Our implementation is based on a discretization of the Riemannian energy functional using B-splines. The boundary value problem for geodesics is solved by a gradient descent on the set of discretized paths and the initial value problem by discrete geodesic calculus  \cite{Rumpf2014}. Our approach is general in that it involves no restriction on the parameters of the metric and allows to factor out rigid transformations. 

We tested our implementation on a dataset of shapes representing various groups of similar physical objects \cite{kimia}. We were able to reconstruct the groups using agglomerative clustering based on the matrix of pairwise geodesic distances between the shapes. We also studied within-group shape variations by performing nonlinear principal component analysis using the Karcher means of the groups as base points. We visualized the results by solving the geodesic equation forward and backward in time along the principal directions. 

Our algorithms performed well on the shapes in the database. Some care was needed to avoid singularities in the spline representations of the contours extracted from the binary images and of the paths initializing the gradient descent of the energy functional. 

In future work, our framework could be applied to other spaces of mappings like manifold-valued curves, embedded surfaces, or more general spaces of immersions (see \cite{Bauer2011a, Bauer2011b} for details and \cite{Bauer2014} for a general overview). A rigorous convergence analysis of the proposed discretization remains to be done.

\section{Mathematical background}\label{sec:background}

We extend the exposition of \cite{BBHM2015} to unparametrized curves. The space of smooth, regular, parametrized curves is defined as
\begin{equation*}
\on{Imm}(S^1,\mathbb R^d)=\left\{c\in C^{\infty}(S^1,\mathbb R^d)\colon \forall \th \in S^1, c_\th(\th) \neq 0 \right\}\,,
\end{equation*}
where $\on{Imm}$ stands for \emph{immersion}. $\on{Imm}(S^1,\R^d)$ is an open subset of the Fr\'echet space $C^\infty(S^1,\R^d)$ and as such itself a Fr\'echet manifold. The tangent space $T_c \on{Imm}(S^1,\R^d)$ at any curve $c$ is the vector space $C^\infty(S^1,\R^d)$. Geometrically, tangent vectors at $c$ are $\R^d$-valued vector fields along $c$.

We denote the Euclidean inner product on $\mathbb R^d$ by $\langle\cdot,\cdot\rangle$. Moreover, for any fixed curve $c$, we denote differentiation and integration with respect to arc length by $D_s=\frac{1}{|c_\theta|}\partial_{\theta}$ and $ds=|c_\theta|d\theta$, respectively.

\begin{definition}
\label{def:sobolev_metric}
Second order Sobolev metrics on $\on{Imm}(S^1,\mathbb R^d)$ are Riemannian metrics of the form
\begin{equation*}
G_c(h,k) = \int_{S^1} a_0\langle h,k \rangle+a_1\langle D_s h,D_s k \rangle+a_2\langle D_s^2 h,D_s^2 k \rangle \,ds \,,
\end{equation*}
where $a_0, a_2 > 0$, $a_1 \geq 0$ are constants and $h,k \in T_c\on{Imm}(S^1,\mathbb R^d)$ tangent vectors.
\end{definition}

\begin{remark}[Choice of the constants]
In the definition of Sobolev metrics, there is freedom in the choice of the relative weighting of the $L^2$-, $H^1$- and $H^2$-parts. This choice influences the geometry of the space of curves and should, in each application, be informed by the data at hand.
\end{remark}

The reparametrization group is the diffeomorphism group of the circle,
\begin{equation*}
\on{Diff}(S^1) = \left\{ \ph \in C^\infty(S^1,S^1) \,:\, \ph \text{ bijective} \right\}\,,
\end{equation*}
which is an infinite-dimensional regular Fr\'echet Lie group \cite{Michor1997}. Reparametrizations act on curves by composition from the right, $(c, \ph) \mapsto c \circ \ph$. The shape space of unparametrized curves is the orbit space $\mc S(S^1,\R^d) = \on{Imm}(S^1,\R^d) / \on{Diff}(S^1)$ of this group action.\footnote{More precisely, we define $\mc S(S^1,\R^d) = \on{Imm}_{\on{f}}(S^1,\R^d) / \on{Diff}(S^1)$, where $\on{Imm}_{\on{f}}(S^1,\R^d)$ is the subset of free immersions, i.e., those upon which $\on{Diff}(S^1)$ acts freely. This is an open and dense subset of $\on{Imm}(S^1,\R^d)$. While important for theoretical reasons, this restriction has no influence on the practical applications of Sobolev metrics.} 

\begin{theorem}[\cite{Michor1991}, Sect. 1.5]\label{thm:shape_space} 
The space $\mc S(S^1,\R^d)$ is a Fr\'echet manifold and the base space of the principal fibre bundle
\begin{equation*}
\pi: \on{Imm}(S^1,\R^d) \to \mc S(S^1,\R^d)\,,\quad c \mapsto c \circ \on{Diff}(S^1)\,,
\end{equation*}
with structure group $\on{Diff}(S^1)$. A Sobolev metric $G$ on $\on{Imm}(S^1,\R^d)$ induces a metric on $\mc S(S^1,\R^d)$, such that the projection $\pi$ is a Riemannian submersion.
\end{theorem}

\begin{remark}[Invariance of the metric]
Thm.~\ref{thm:shape_space} hinges on the invariance of Sobolev metrics with respect to reparametrizations. Sobolev metrics are also invariant with respect to translations and rotations, but in general not to scalings. The lack of scale invariance can be addressed by introducing weights depending on the length $\ell_c$ of the curve $c$ as follows:
\begin{equation*}
\widetilde G_c(h,k) = \int_{S^1} \frac{a_0}{\ell^3_c}\langle h,k \rangle+\frac{a_1}{\ell_c} \langle D_s h,D_s k \rangle+a_2\ell_c \langle D_s^2 h,D_s^2 k \rangle\, ds\,.
\end{equation*}
\end{remark}

Deformations of curves are smooth paths $c\colon [0,1]\to\on{Imm}(S^1,\mathbb R^d)$. Their velocity is $c_t$, the subscript $t$ denoting differentiation. The length of a path $c$ is
\begin{equation*}
L(c) = \int_0^1 \sqrt{G_{c(t)}(c_t(t),c_t(t))} \, dt\,.
\end{equation*}
The distance between two shapes $\pi(c_0)$ and $\pi(c_1)$ in $\mc S(S^1,\R^d)$ with respect to the metric induced by $G$ is the infimum over all paths between $c_0$ and the orbit $c_1 \circ \on{Diff}(S^1)$, i.e.,
\begin{equation*}
\on{dist}(c_0, c_1) = \inf_{\substack{c(0) = c_0 \\ c(1) \in c_1 \circ \on{Diff}(S^1)}} L(c)\,.	
\end{equation*}
Geodesics on $\on{Imm}(S^1,\R^d)$ are critical points of the energy functional
\begin{equation*}
E(c) = \frac12\int_0^1 G_{c(t)}\big(c_t(t),c_t(t)\big) \, dt\,.
 \label{eq: EnergyFunctional}
\end{equation*}
The spaces $\mc S(S^1,\R^d)$ and $\on{Imm}(S^1,\R^d)$ are related by a Riemannian submersion and geodesics on $\mc S(S^1,\R^d)$ can be lifted to horizontal geodesics on $\on{Imm}(S^1,\R^d)$. Conversely, the projection of a length-minimizing path between $c_0$ and $c_1 \circ \on{Diff}(S^1,\R^d)$ is a geodesic on $\mc S(S^1,\R^d)$.

The space $\mc S(S^1,\R^d)$ equipped with a Sobolev metric possesses some nice completeness properties, which are summarized in the following theorem.

\begin{theorem}[Minimizing geodesics,  \cite{Bruveris2014,Bruveris2014b_preprint}]
Let $a_0, a_2 > 0$ and $a_1 \geq 0$. Then, given two curves $c_0, c_1$ in the same connected component of $\on{Imm}(S^1,\R^d)$, there exists a minimizing geodesic connecting them. Furthermore, there exists a minimizing geodesic connecting the shapes $\pi(c_0)$ and $\pi(c_1)$ in $\mc S(S^1,\R^d)$.
\end{theorem}

\begin{remark}[Elastic metrics]
Closely related to the Sobolev metrics described here is the family of elastic metrics \cite{Mio2007}, which in the planar case is given by
\begin{equation*} 
 G_c(h,k) = \int_{S^1} a^2\langle D_s h,n \rangle \langle D_s k,n \rangle+b^2\langle D_s h,v \rangle\langle D_s k,v \rangle ds\,.
\end{equation*}
Here $a,b$ are constants and $v,n$ denote the unit tangent and normal vectors to $c$. Two special cases deserve to be highlighted: for $a=1$, $b=\frac 12$  \cite{Jermyn2011} and $a=b$ \cite{Michor2008a} there exist nonlinear transforms, the square root velocity transform and the basic mapping, that greatly simplify the numerical computation of geodesics. Both of these metrics have been applied to a variety of problems in shape analysis. We note that the elastic metric with $a=b$ corresponds to a first order Sobolev metric as in Def.~\ref{def:sobolev_metric} with $a_0=a_2=0$ and $a_1 = a^2=b^2$. As it has no $L^2$-part, it is a Riemannian metric only on the space of curves modulo translations.
\end{remark}

\section{Numerical implementation}

\subsection{Discretization} 

We discretize curves using B-splines; $c = \sum_{j=1}^{N_\theta} d_j C_j(\theta)$, where $C_j$ are the B-splines of degree $n_{\theta}$, defined on a uniform periodic knot sequence, with all knots of multiplicity one. Observe that $C_j \in C^{n_{\theta}-1}([0,2\pi])$. A path of curves can then be represented using tensor product B-splines, i.e.,
\begin{equation}
c(t,\theta) = \sum_{i=1}^{N_t} \sum_{j=1}^{N_\theta} d_{i,j} B_i(t)C_j(\theta)\,.
\label{eq:TensorProductPath} 
\end{equation}  
Here $B_i$ are the B-splines defined on the interval $[0,1]$, with uniform knots and full multiplicity at the end points, $B_i \in C^{n_t-1}([0,1])$. This implies that the boundary curves are given by $c(0,\theta) = \sum_{j=1}^{N_\theta} d_{1,j} C_j(\theta)$ and $c(1,\theta) = \sum_{j=1}^{N_\theta} d_{N_{t}, j} C_j(\theta)$.

Under the identification of $S^1$ with $\mathbb{R}/[0,2\pi]$, diffeomorphisms $\psi\colon S^1 \to S^1$ can be written as $\psi = \operatorname{id} + \phi$, where  $\phi$ is a periodic function. We choose to discretize $\phi(\theta) = \sum_{i=1}^{N_{\phi}} \phi_i D_i(\theta)$, where $D_i$ are B-splines of degree $n_\phi$, defined on a uniform periodic knot sequence, similarly to $C_j$. The identity can be written in a B-spline basis using the \emph{Greville abscissas} $\xi_i$, i.e., $\operatorname{id} = \sum_{i=1}^{N_\phi} \xi_i D_i$. Due to the positivity of B-splines, the condition $\psi'>0$ ensuring that $\psi$ is a diffeomorphism takes the form 
\begin{equation}
\phi_{i-1} - \phi_i < \xi_{i} - \xi_{i-1}\,.
\label{eq:PhiDiffeoCondition}
\end{equation}  
To speed up convergence, we introduce an additional variable $\alpha \in \mathbb R$ representing constant shifts of the reparametrization. The resulting redundancy is eliminated by the constraint
\begin{equation}\label{eq: shift constraint}
	\sum_{i=1}^{N_\phi} \phi_i = 0\,.
\end{equation}
To compute the energy of paths \eqref{eq:TensorProductPath}, we use Gaussian quadrature on each interval between subsequent knots. Evaluating paths (and their derivatives) at the quadrature points is a simple multiplication of the spline collocation matrix with the vector of control points. 

\subsection{Geodesics and Karcher means}\label{sec:boundary}

From now on we work with plane curves ($d=2$). The boundary value problem for geodesics consists of minimizing the discretized energy \eqref{eq: EnergyFunctional} over all paths $c\colon[0,1] \times [0,2\pi] \to \mathbb{R}^2$, mappings $\psi=\on{id}+\phi\colon S^1\to S^1$, shifts of the reparametrization $\alpha \in \mathbb R$, and rotations $R_\beta$ around the origin by an angle $\beta \in [0,2\pi)$, subject to the constraints \eqref{eq:PhiDiffeoCondition}, \eqref{eq: shift constraint}, and
\begin{equation*}
	c(0, \cdot) = c_0(\cdot), \qquad
	c(1, \cdot) = R_\beta(c_1(\psi(\cdot)-\alpha)+v)\,,
\end{equation*}
where $c_0,c_1\colon[0,2\pi] \to \mathbb{R}^2$ are given boundary curves. This is a finite dimensional constrained optimization problem, which we solve using Matlab's interior point method \texttt{fmincon}. We achieved major performance improvements by fine-tuning the implementations of the gradient and hessian of the energy functional.

In the initial value problem for geodesics, the initial value and initial velocity of the sought geodesic are given. As described in Sect.~\ref{sec:background}, geodesics on $\mc S(S^1,\R^d)$ can be lifted to horizontal geodesics on $\on{Imm}(S^1,\R^d)$. Thus, solving the geodesic initial value problem for unparametrized curves reduces to solving the problem for parametrized curves with horizontal initial velocities. To solve the latter problem, we use the time-discrete variational geodesic calculus of \cite{Rumpf2014} as described in \cite{BBHM2015}.

The Karcher mean $\overline{c}$ of a set $\{c_1,\dots,c_n\}$ of curves is defined as the minimizer of
\begin{equation}\label{eq: karcher functional}
F(c) = \frac{1}{n} \sum_{j=1}^n \on{dist}(c, c_j)^2\,.
\end{equation}
It can be computed by iteratively solving initial and boundary value problems for geodesics. We refer to \cite{Pennec2006b} for more details.

\section{Numerical examples}

\subsection{Data acquisition and setup}

\begin{figure}
\centering
\includegraphics[width=0.1\textwidth,angle=270]{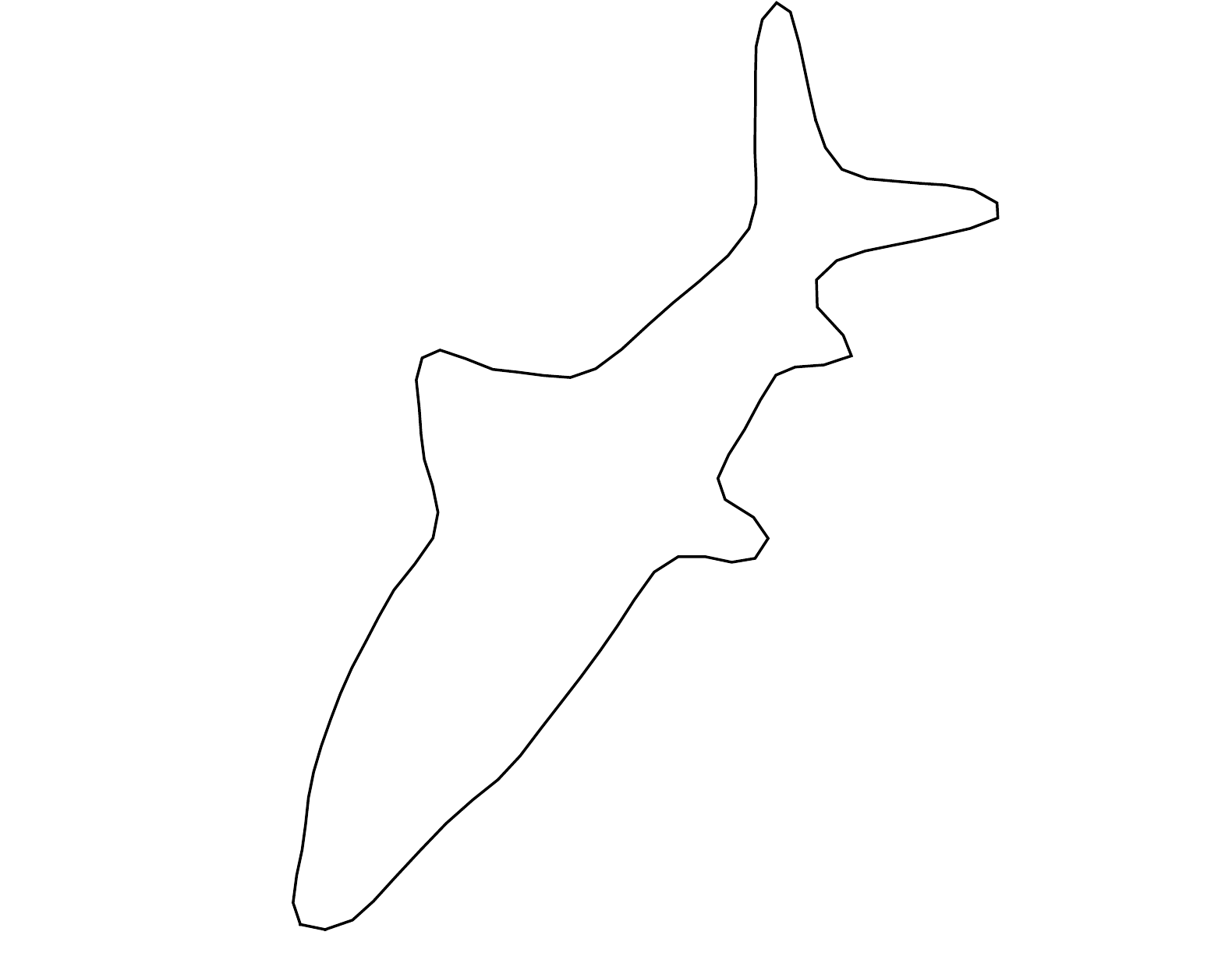}
\includegraphics[width=0.1\textwidth,angle=270]{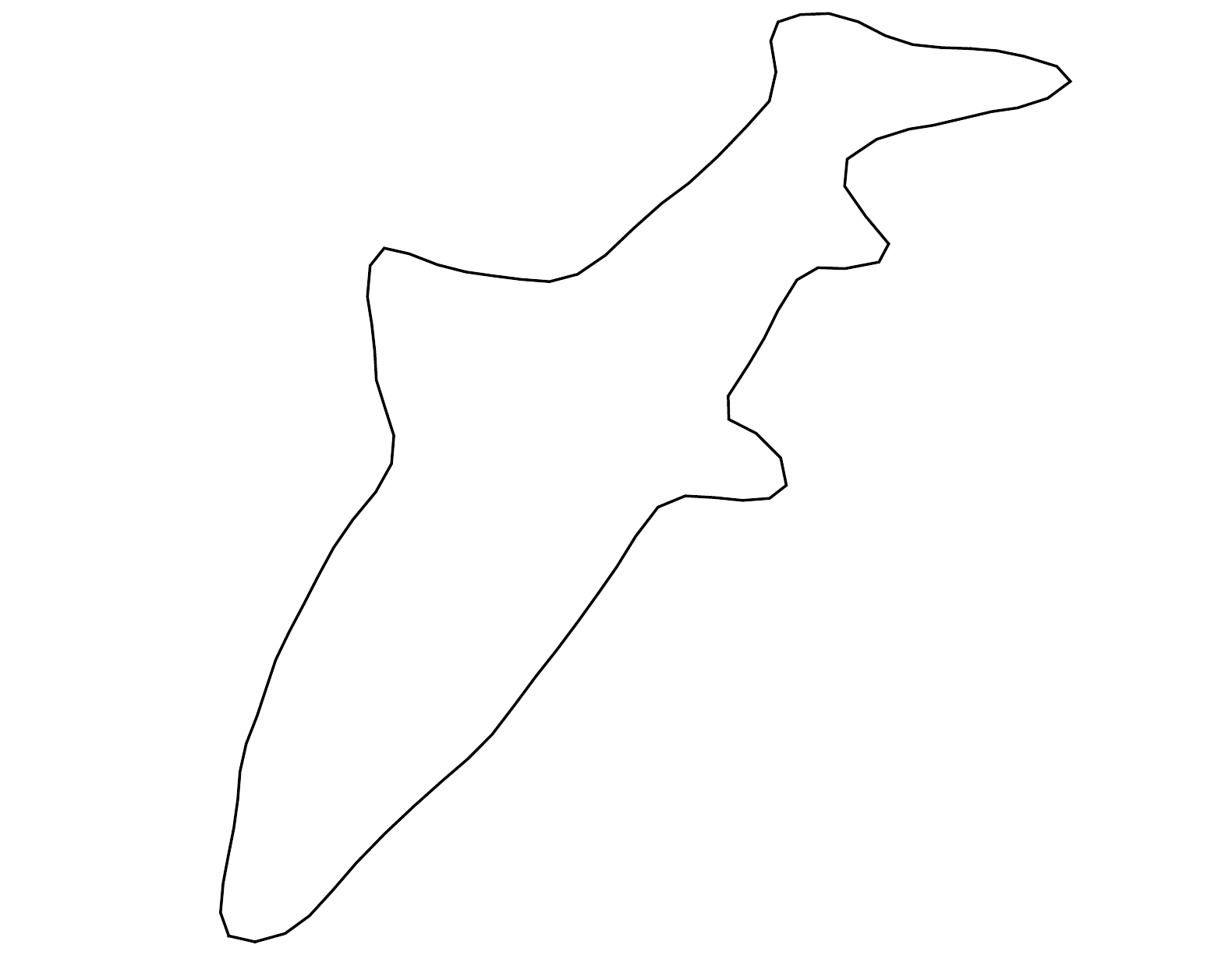}
\includegraphics[width=0.1\textwidth,angle=270]{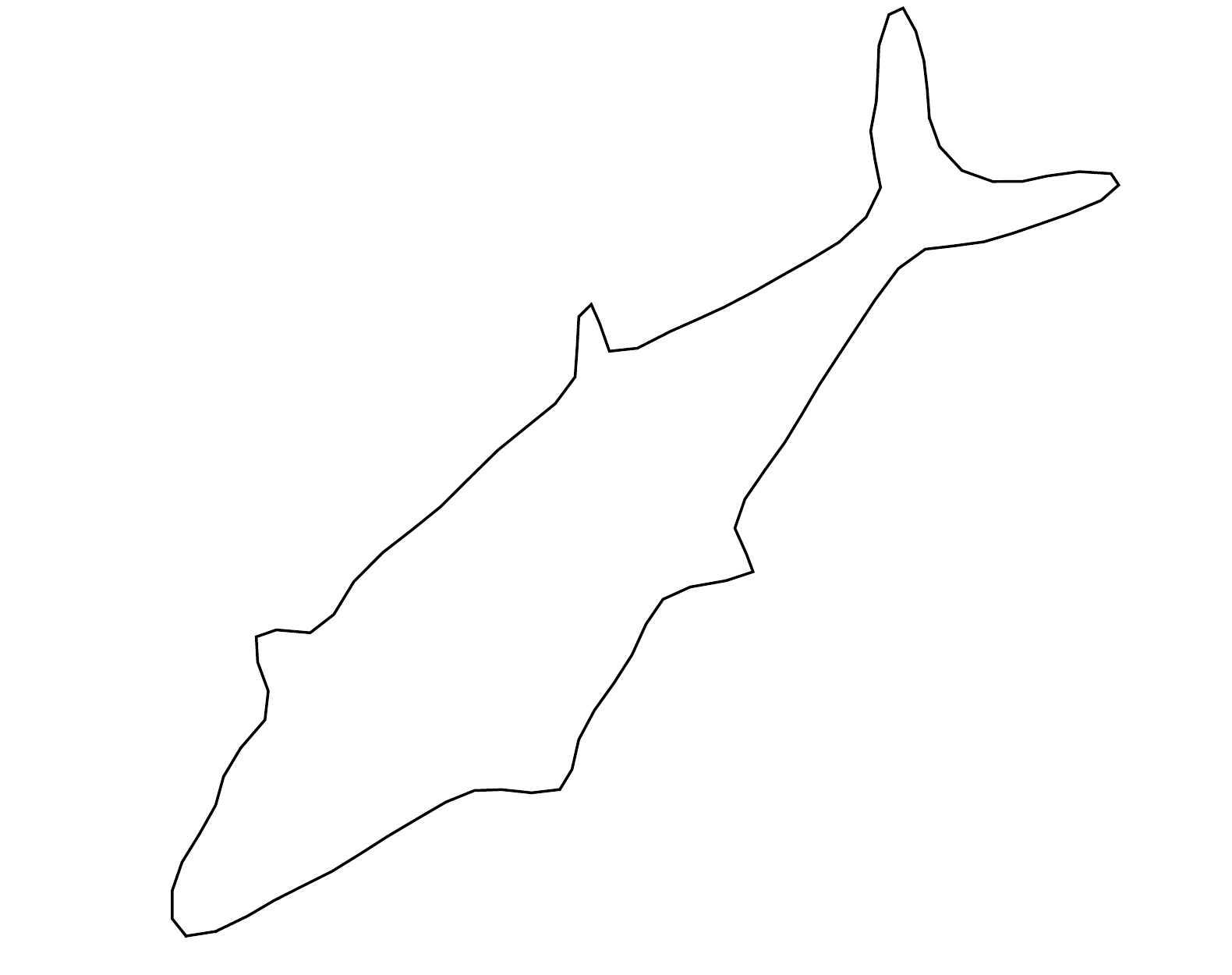}
\includegraphics[width=0.1\textwidth,angle=270]{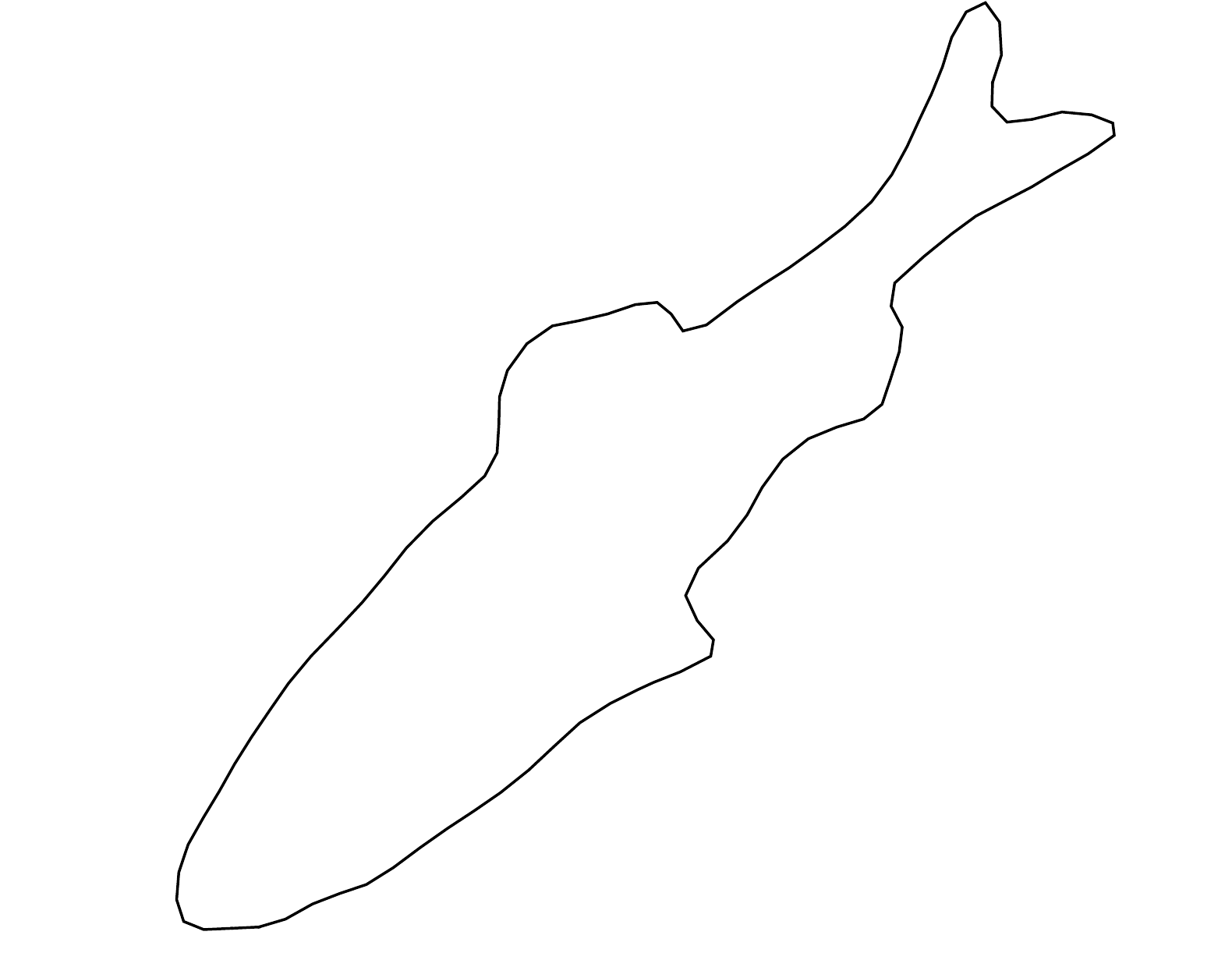}
\includegraphics[width=0.1\textwidth,angle=270]{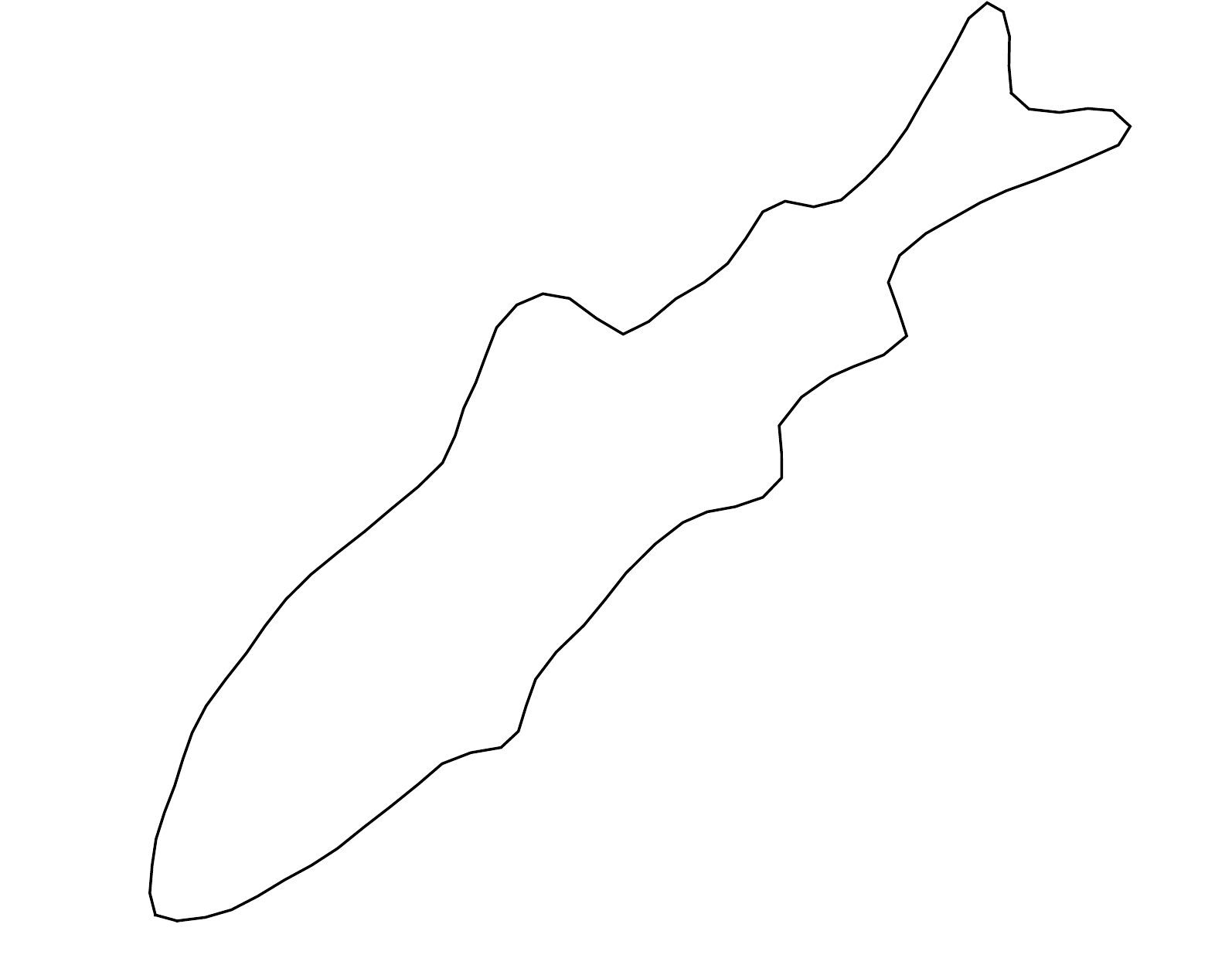}
\includegraphics[width=0.1\textwidth,angle=270]{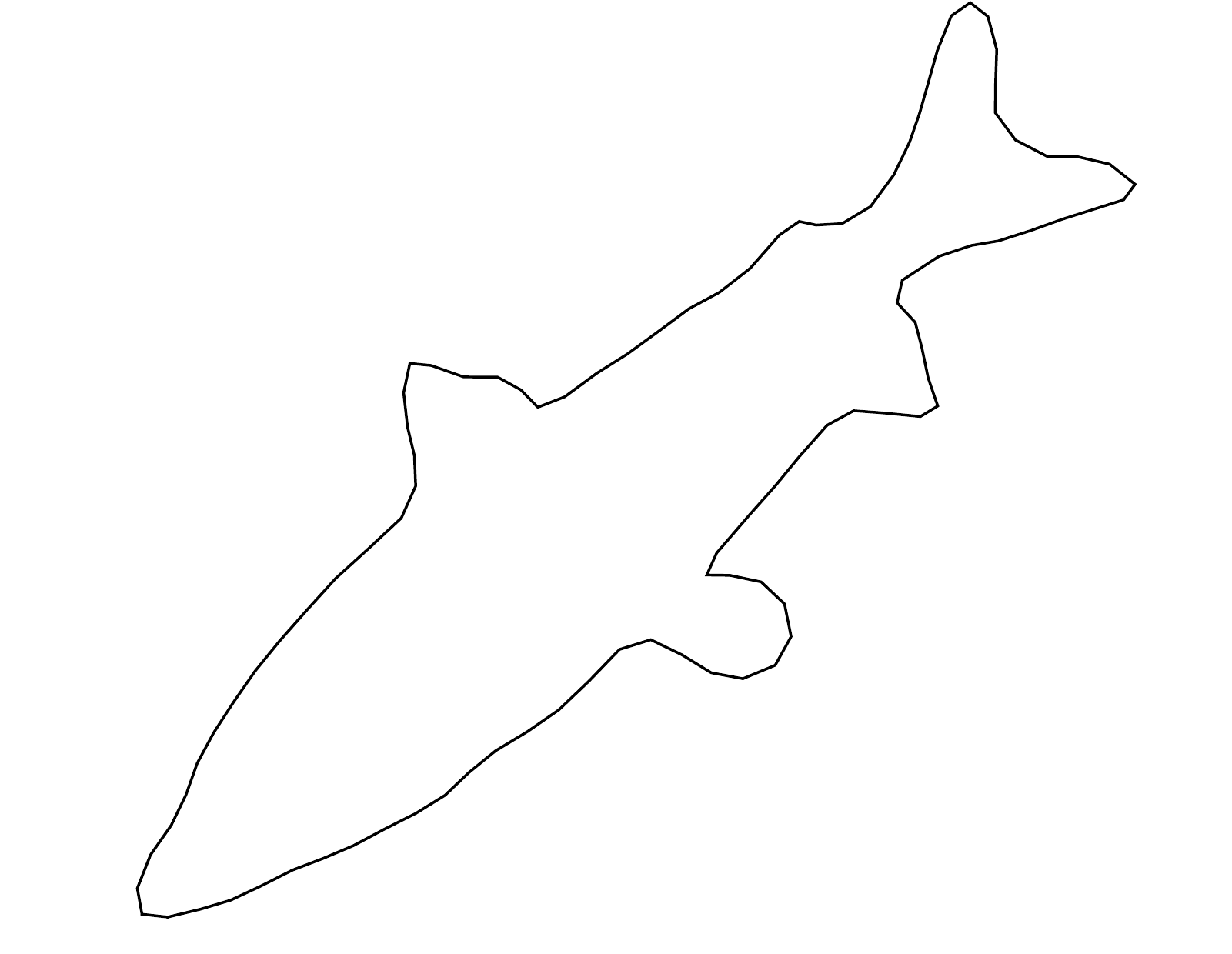}
\includegraphics[width=0.1\textwidth,angle=270]{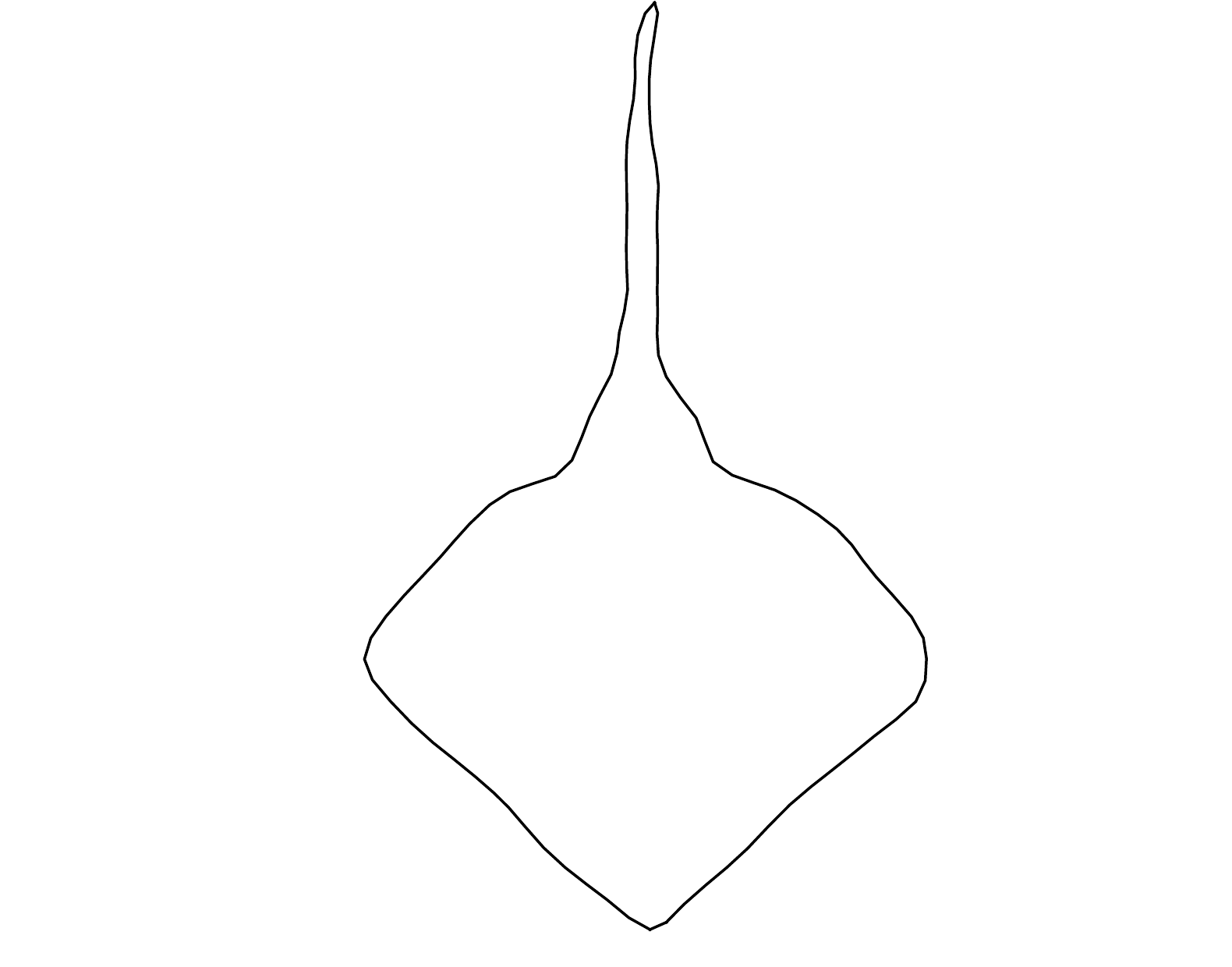}
\includegraphics[width=0.1\textwidth,angle=270]{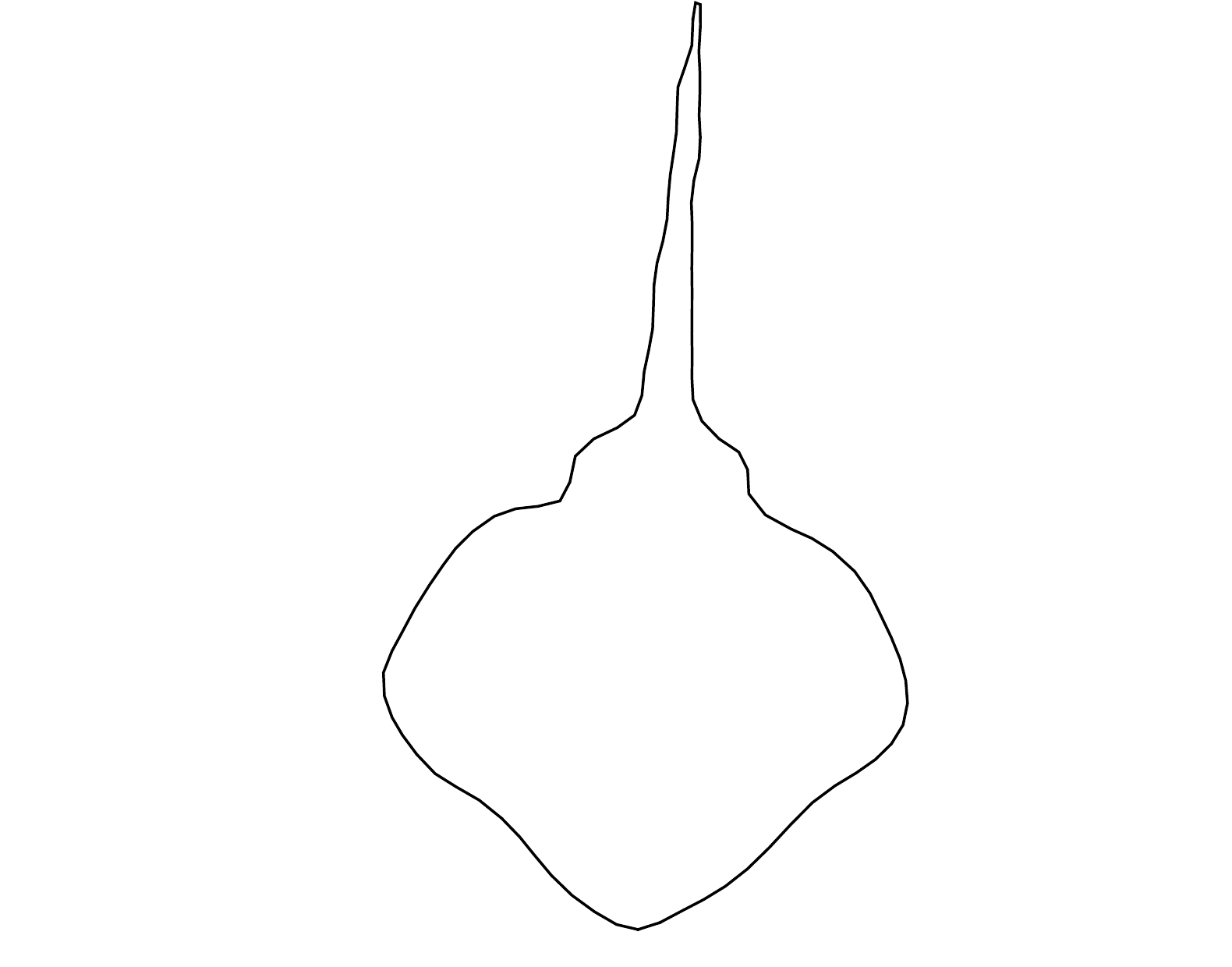}
\includegraphics[width=0.1\textwidth,angle=270]{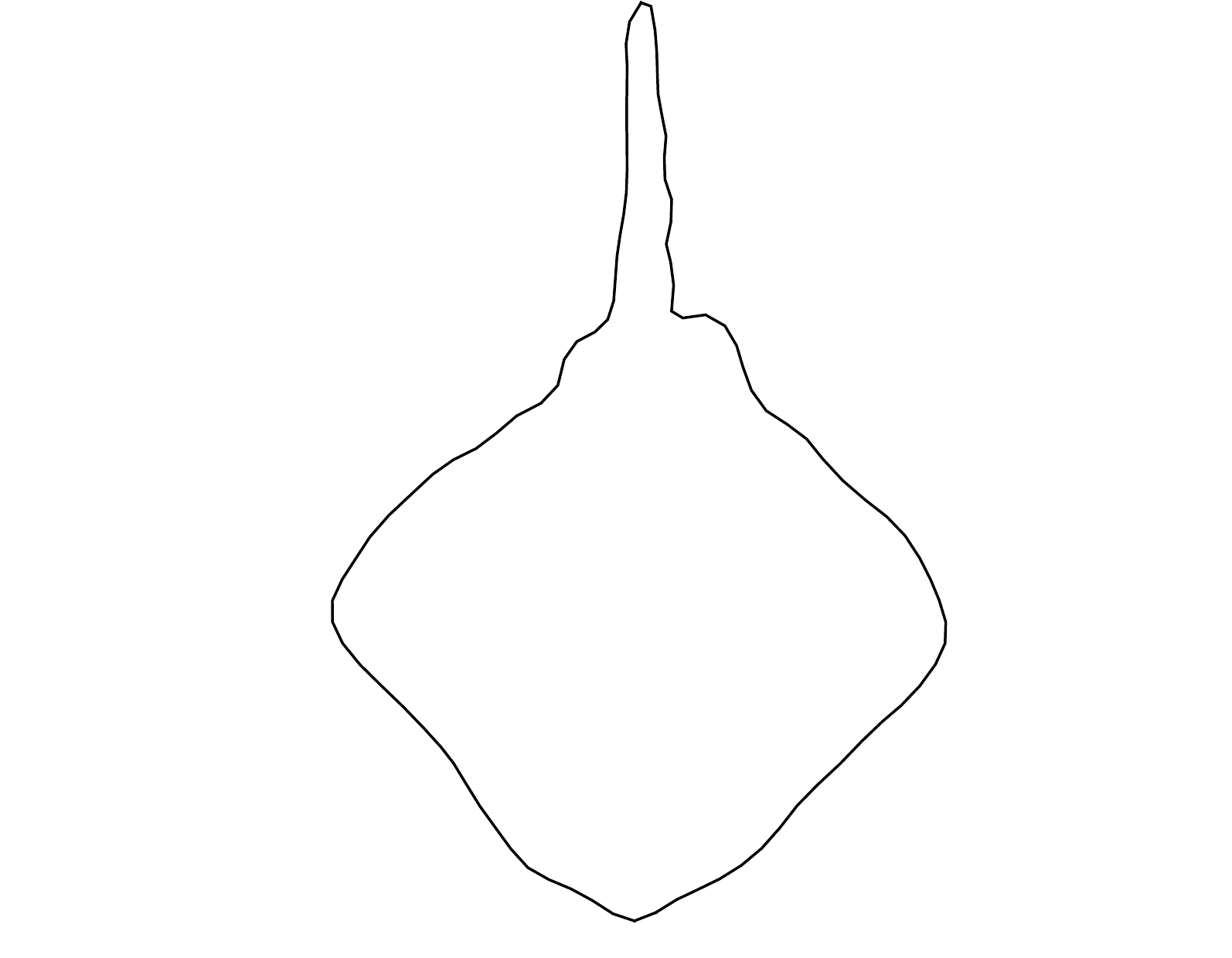}
\includegraphics[width=0.1\textwidth,angle=270]{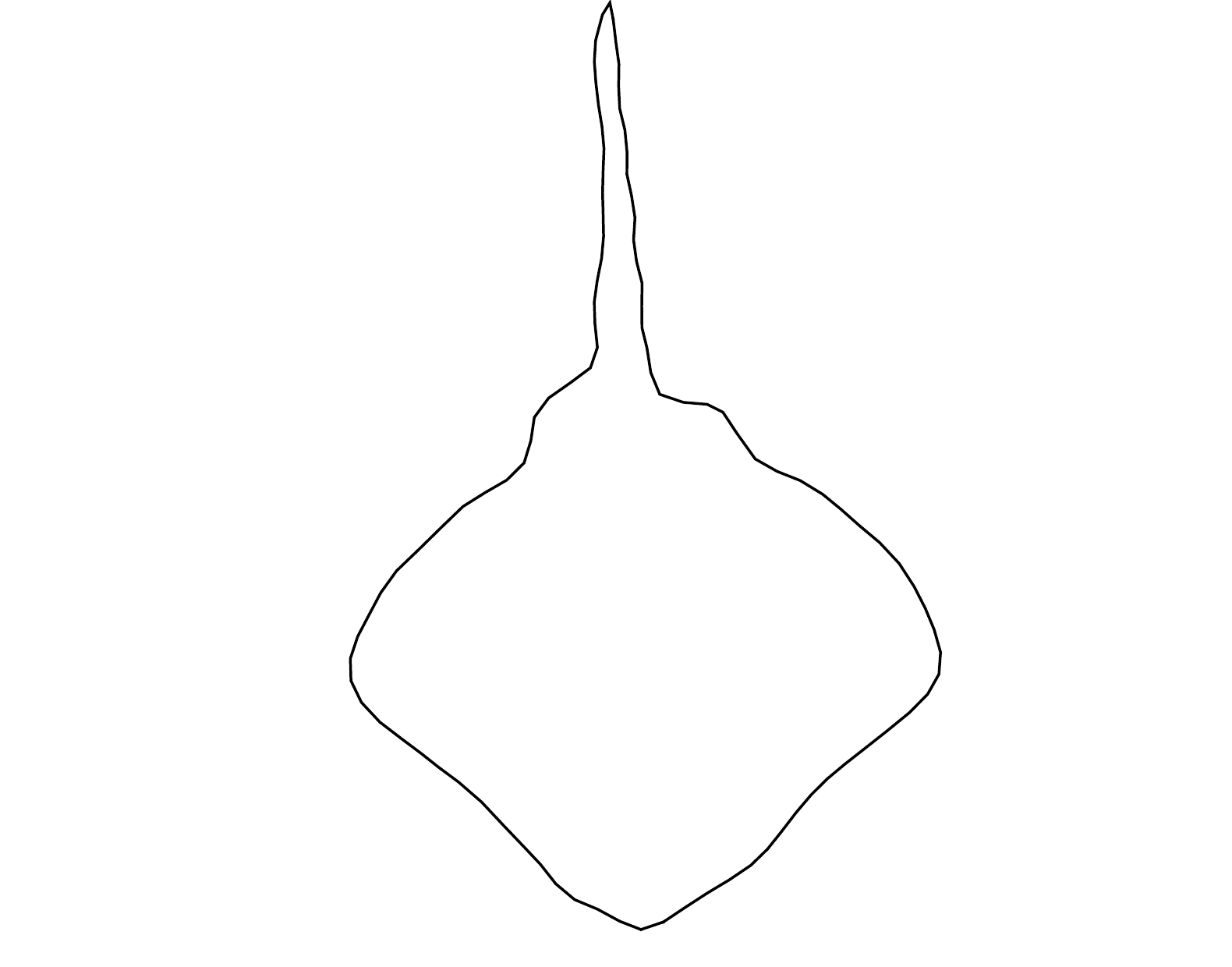}
\includegraphics[width=0.1\textwidth,angle=270]{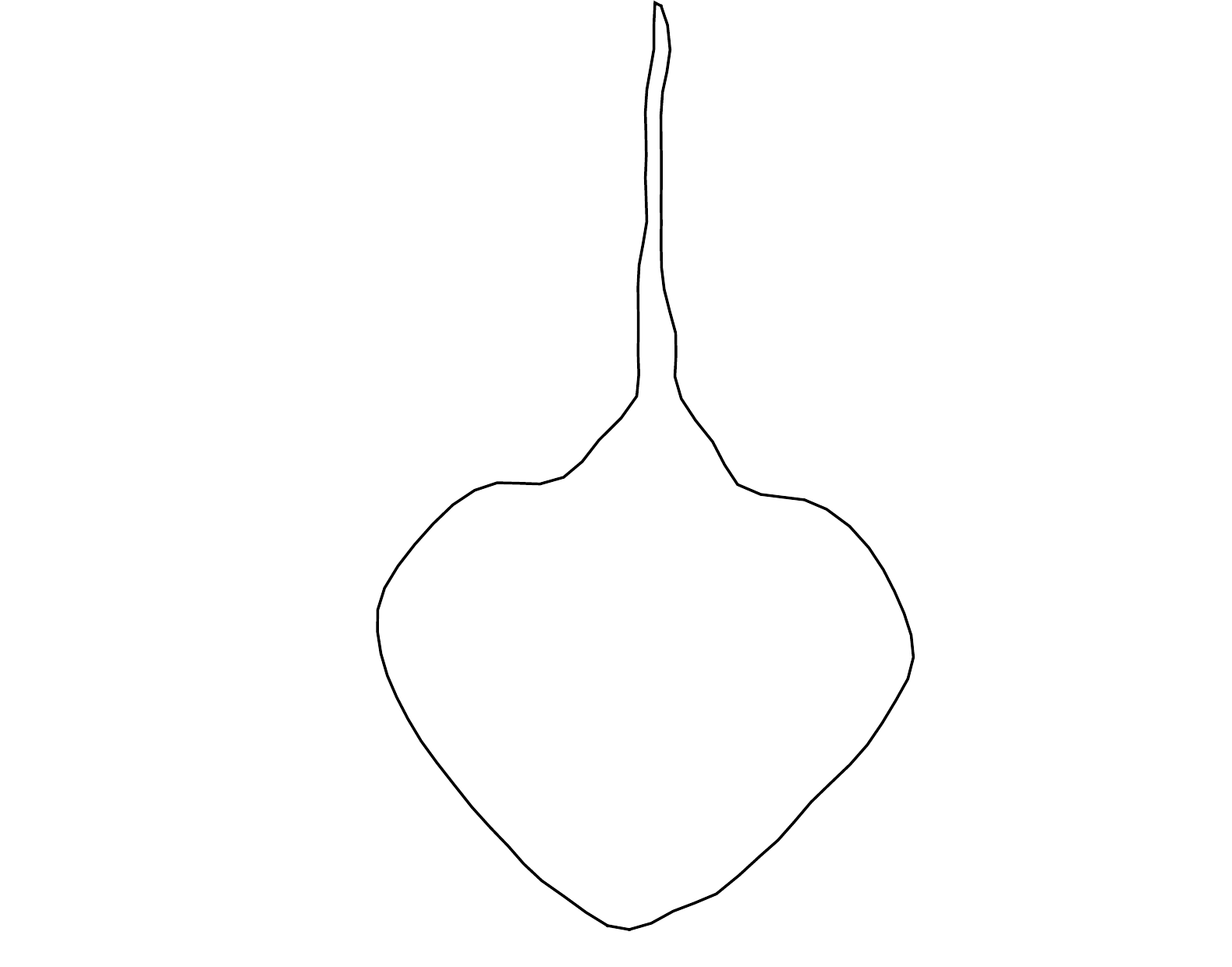}
\includegraphics[width=0.1\textwidth,angle=270]{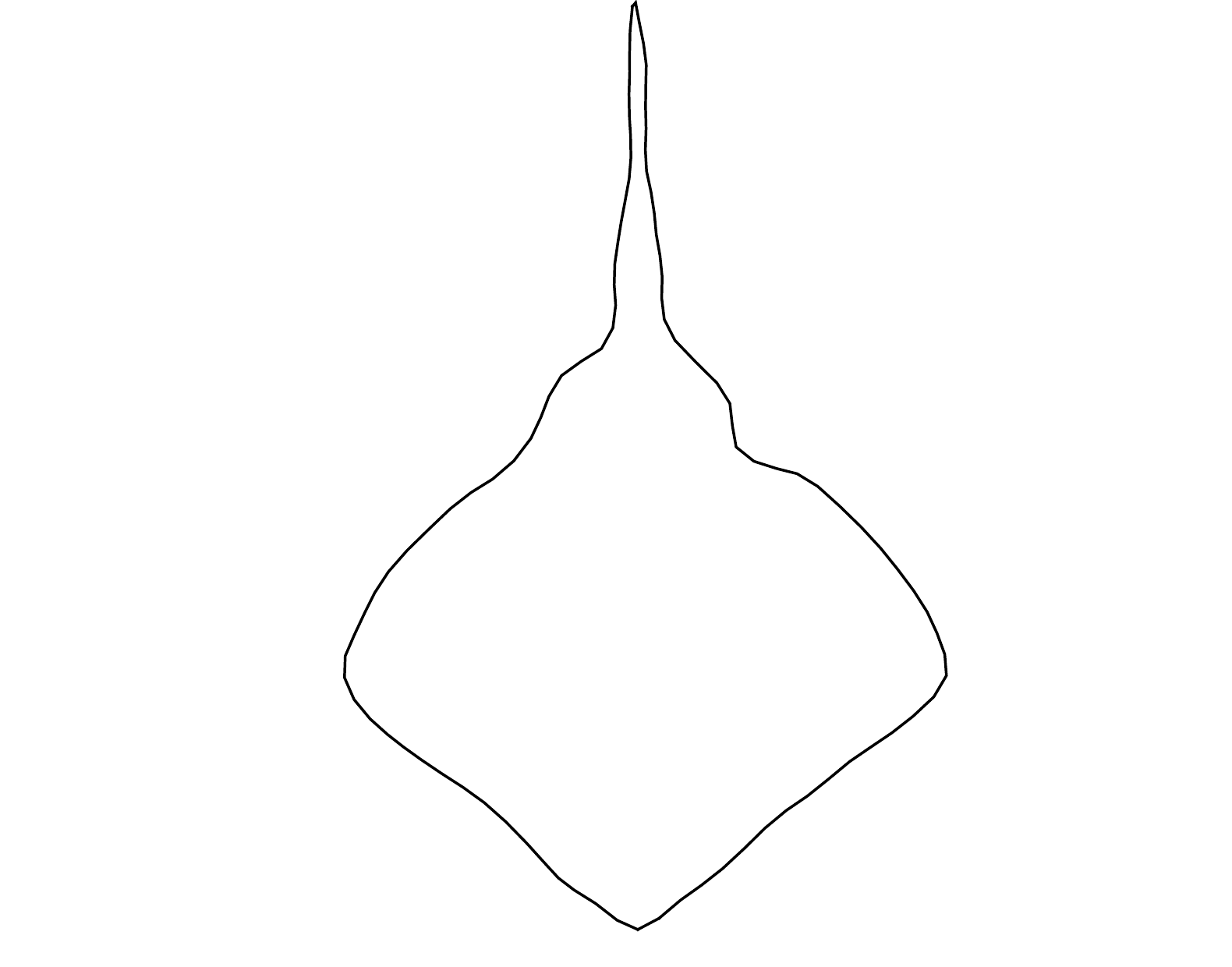}
\includegraphics[width=0.1\textwidth,angle=270]{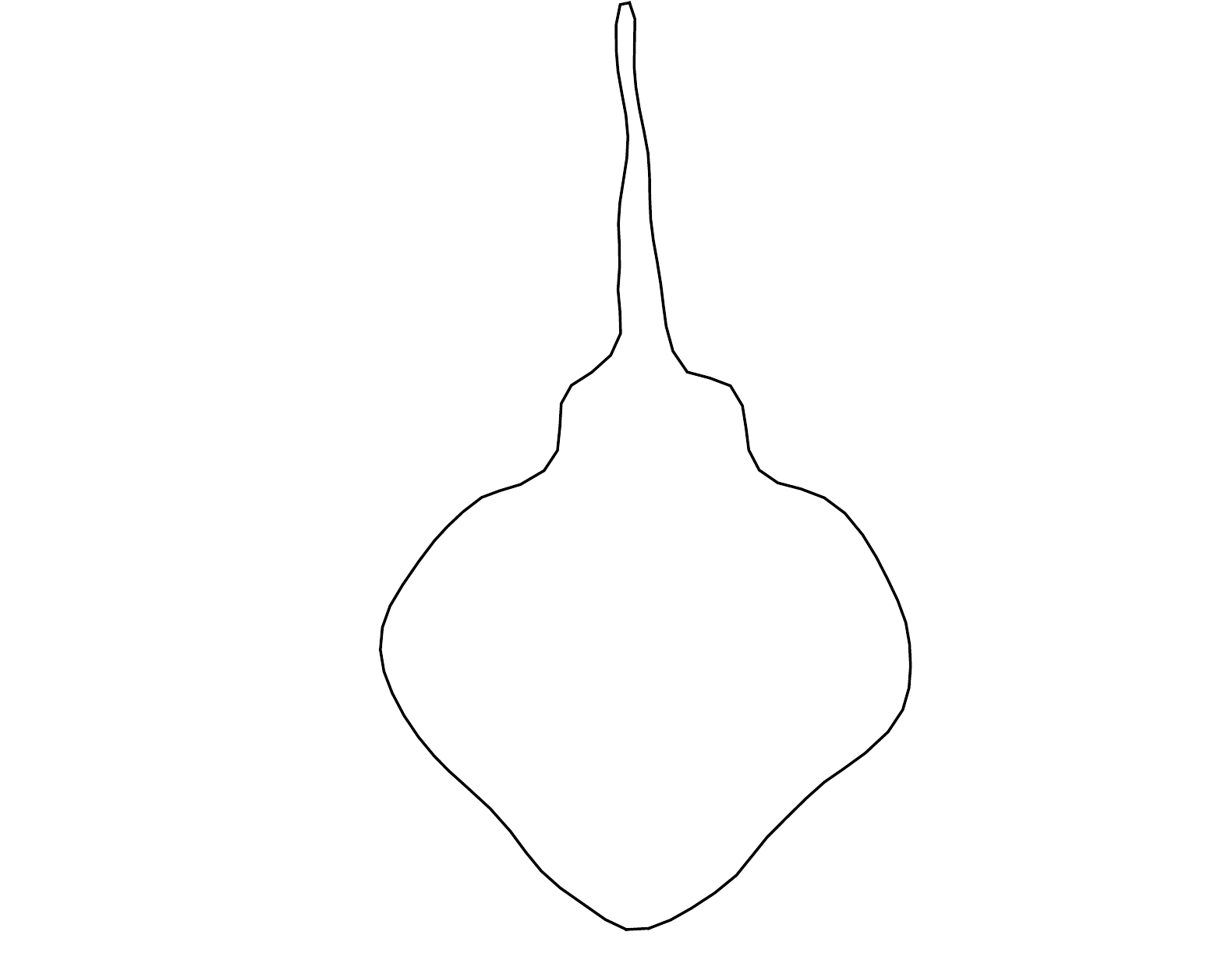}
\includegraphics[width=0.1\textwidth,angle=270]{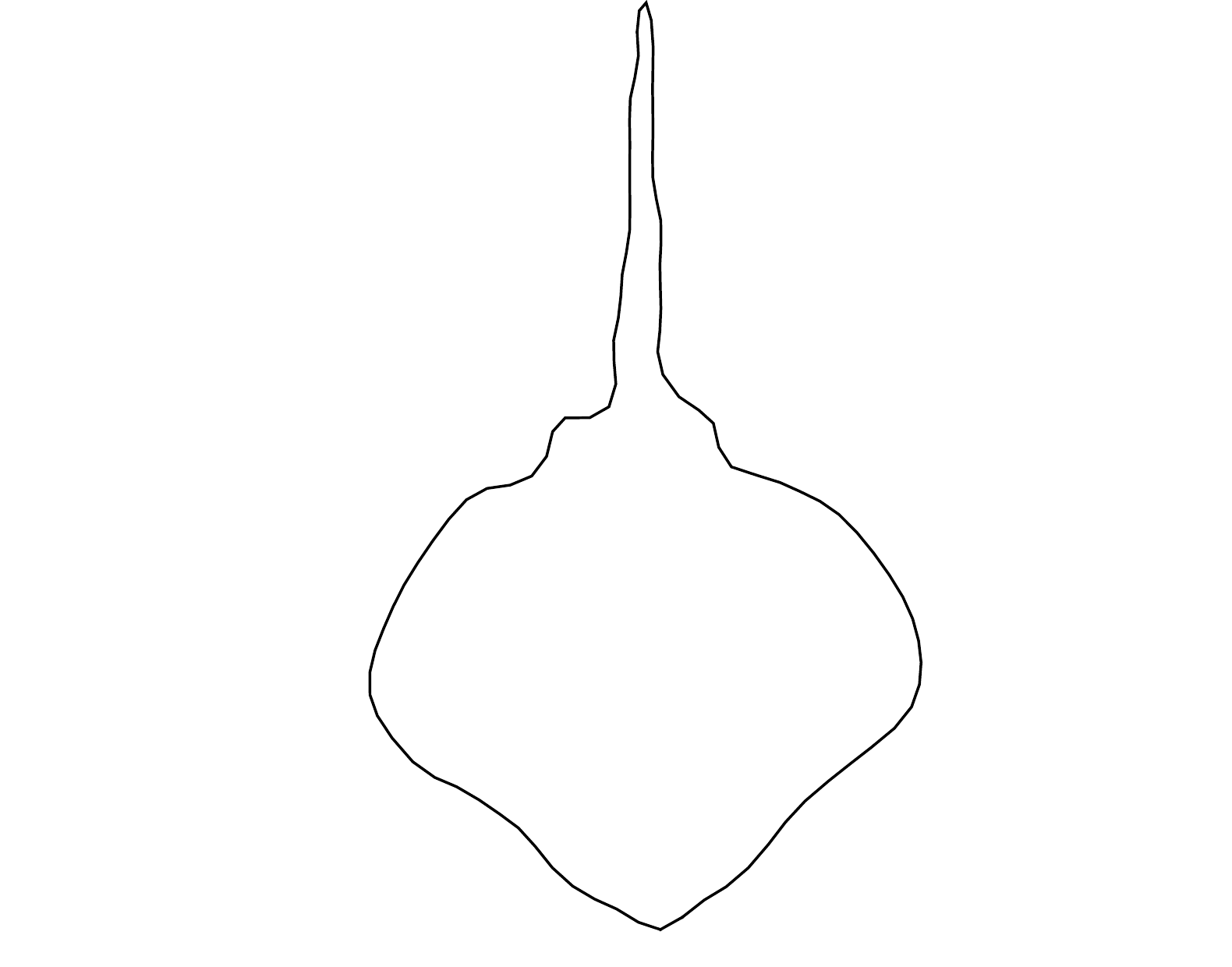}
\includegraphics[width=0.1\textwidth,angle=270]{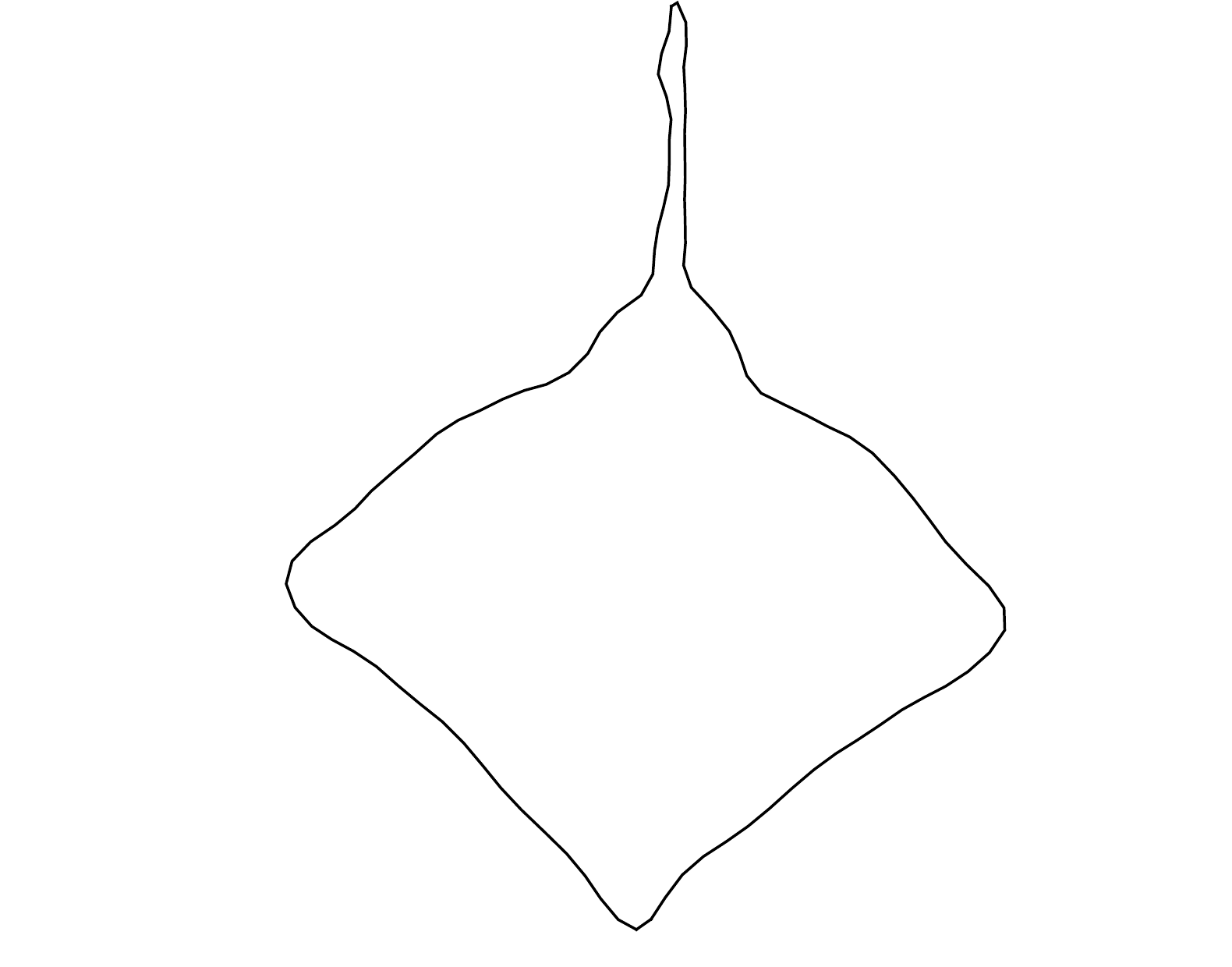}
\includegraphics[width=0.1\textwidth,angle=270]{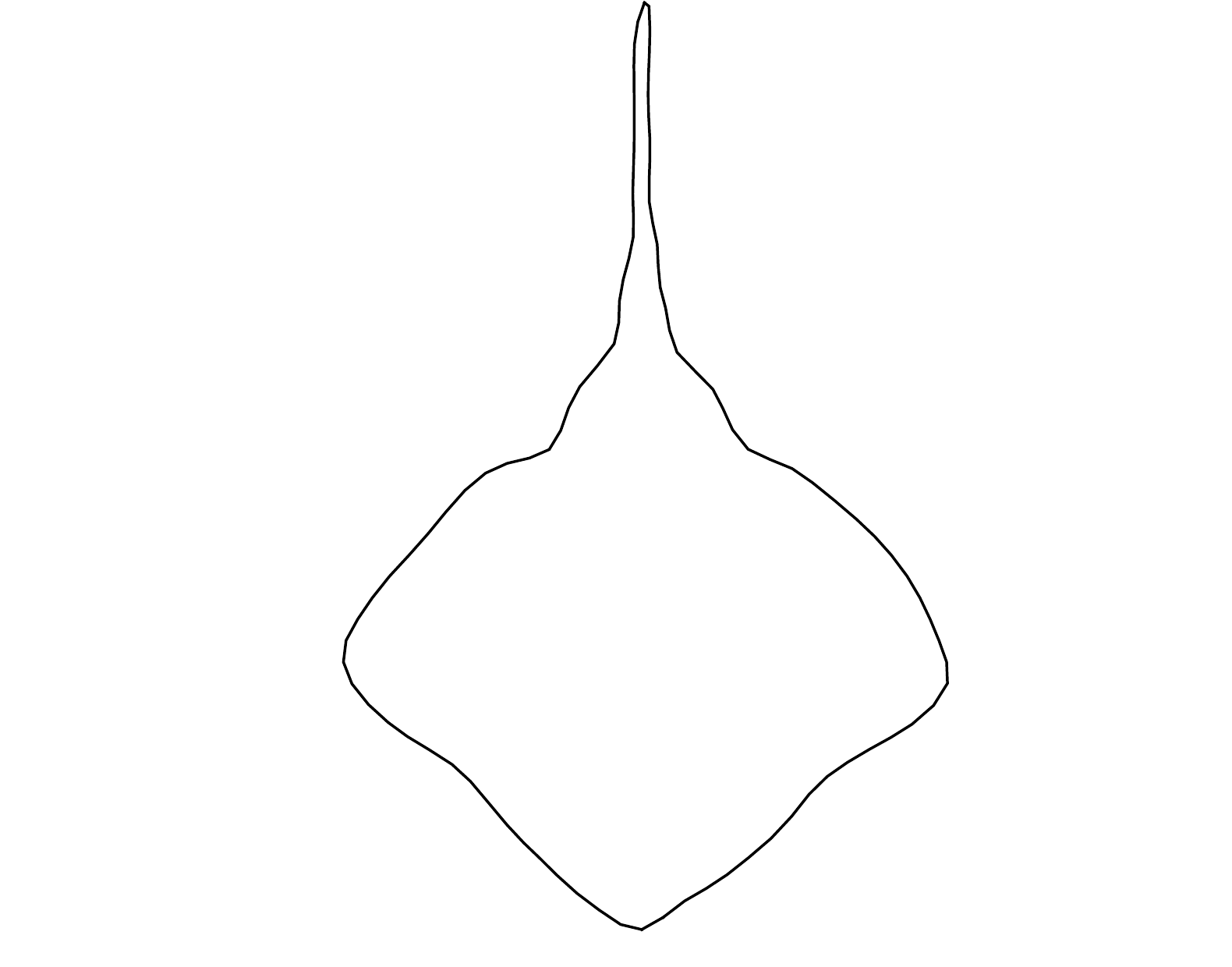}
\includegraphics[width=0.1\textwidth,angle=270]{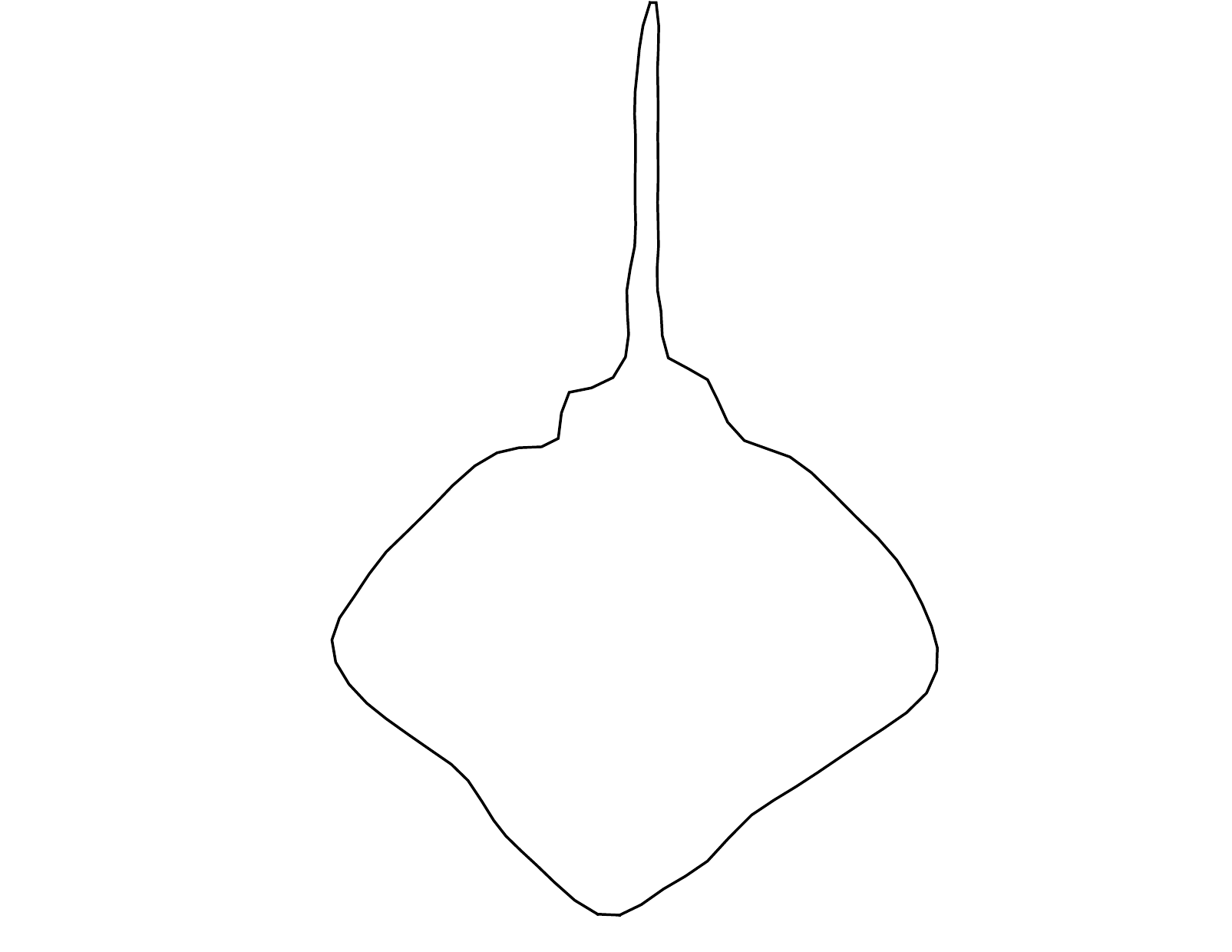}
\includegraphics[width=0.1\textwidth,angle=270]{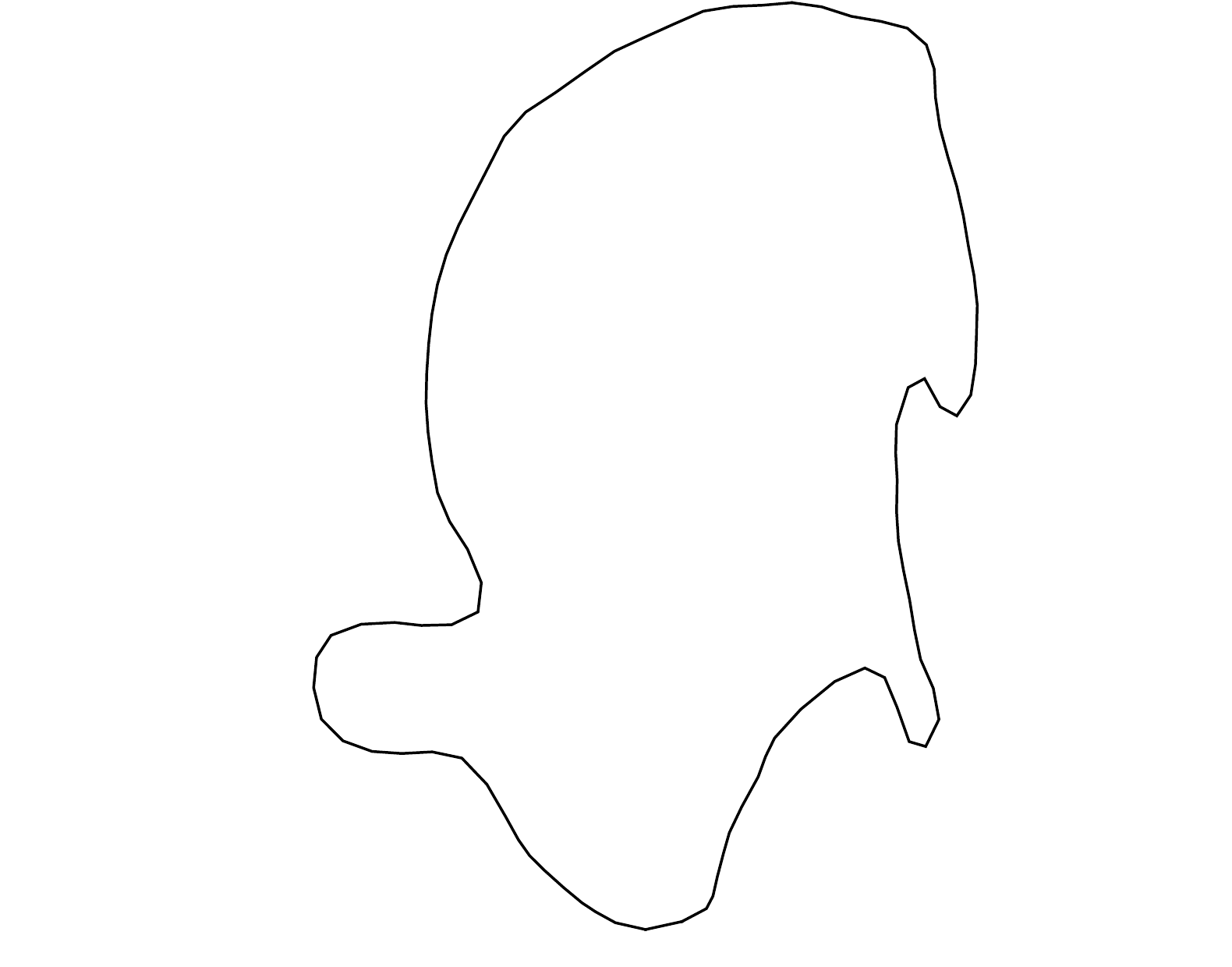}
\includegraphics[width=0.1\textwidth,angle=270]{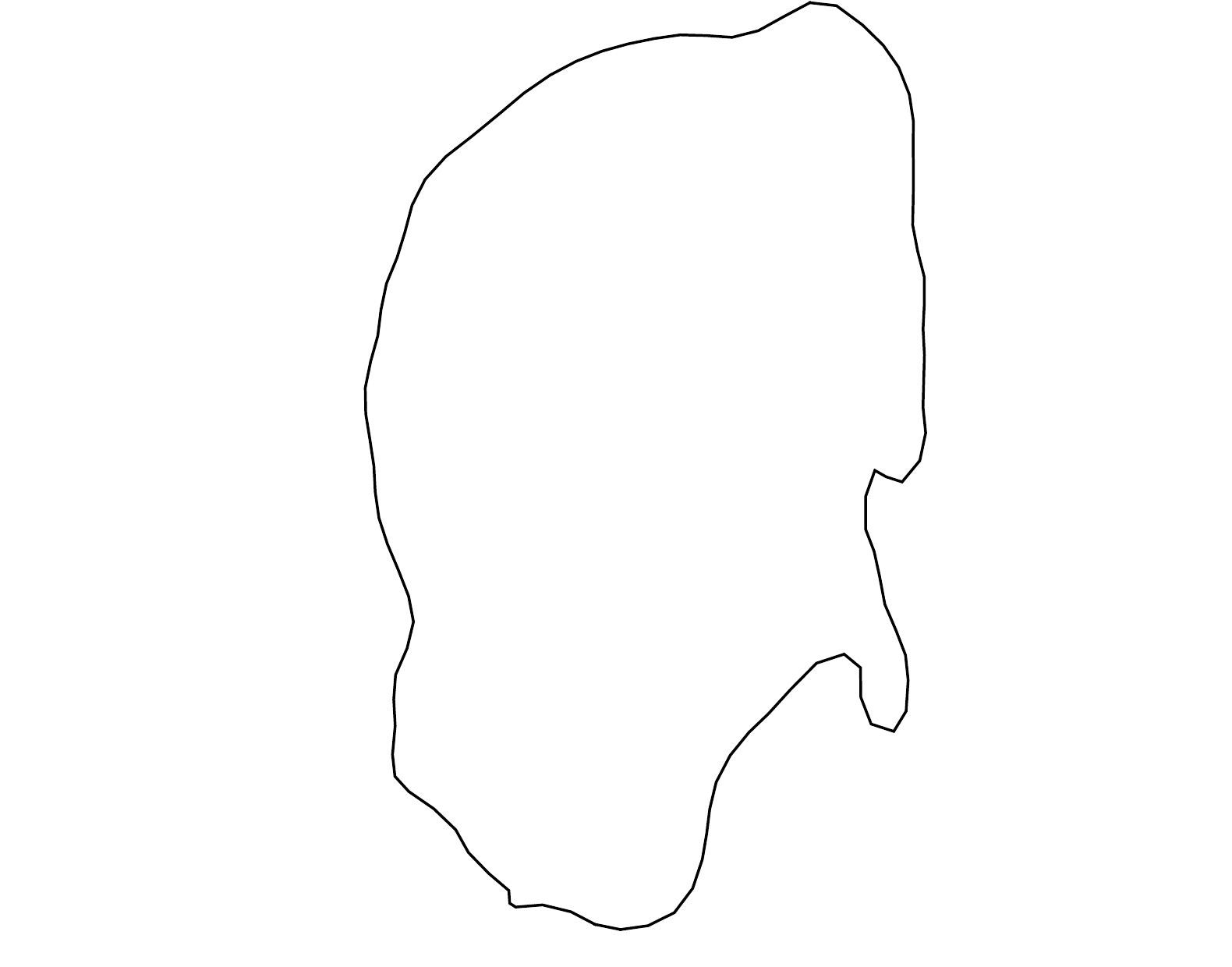}
\includegraphics[width=0.1\textwidth,angle=270]{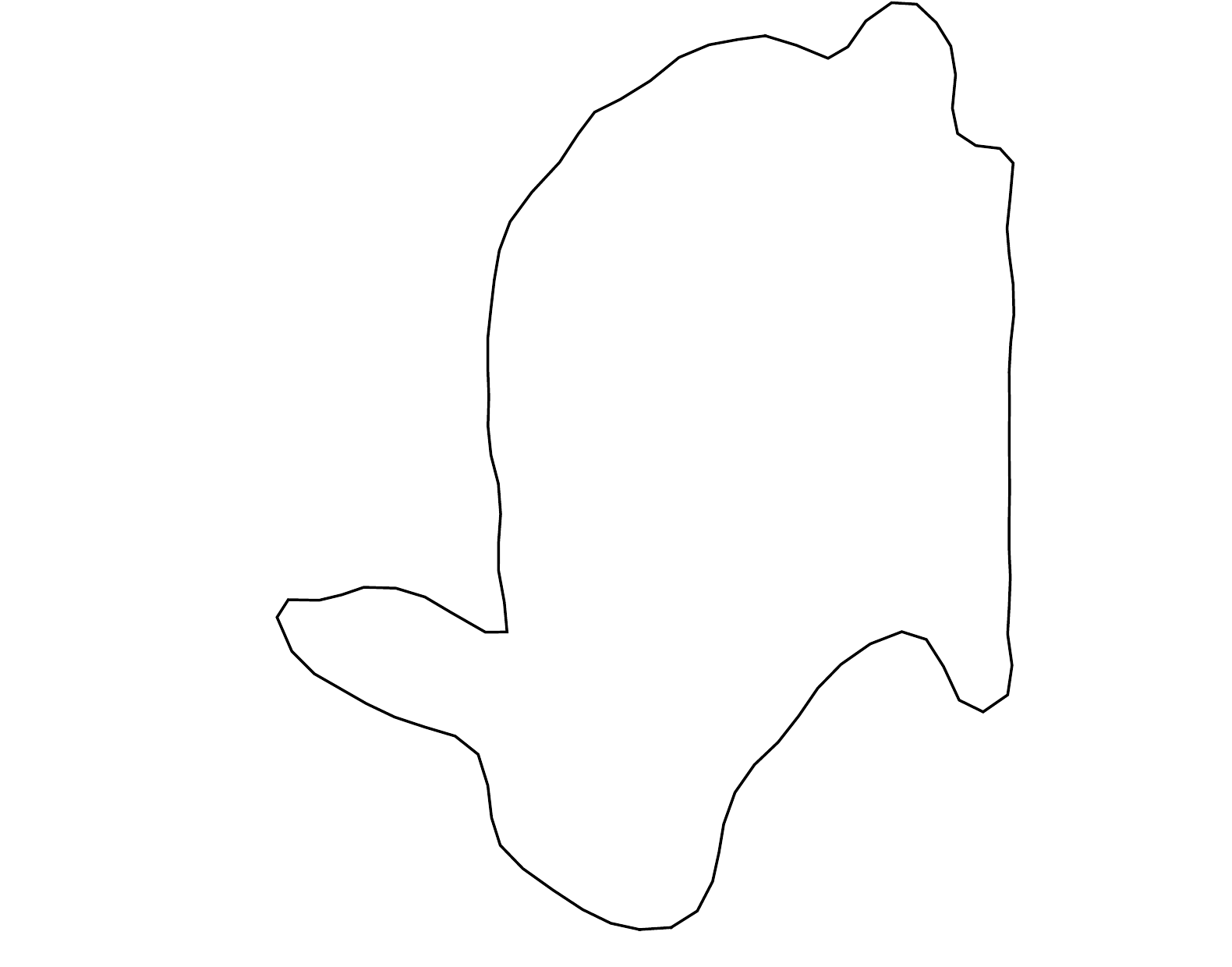}
\includegraphics[width=0.1\textwidth,angle=270]{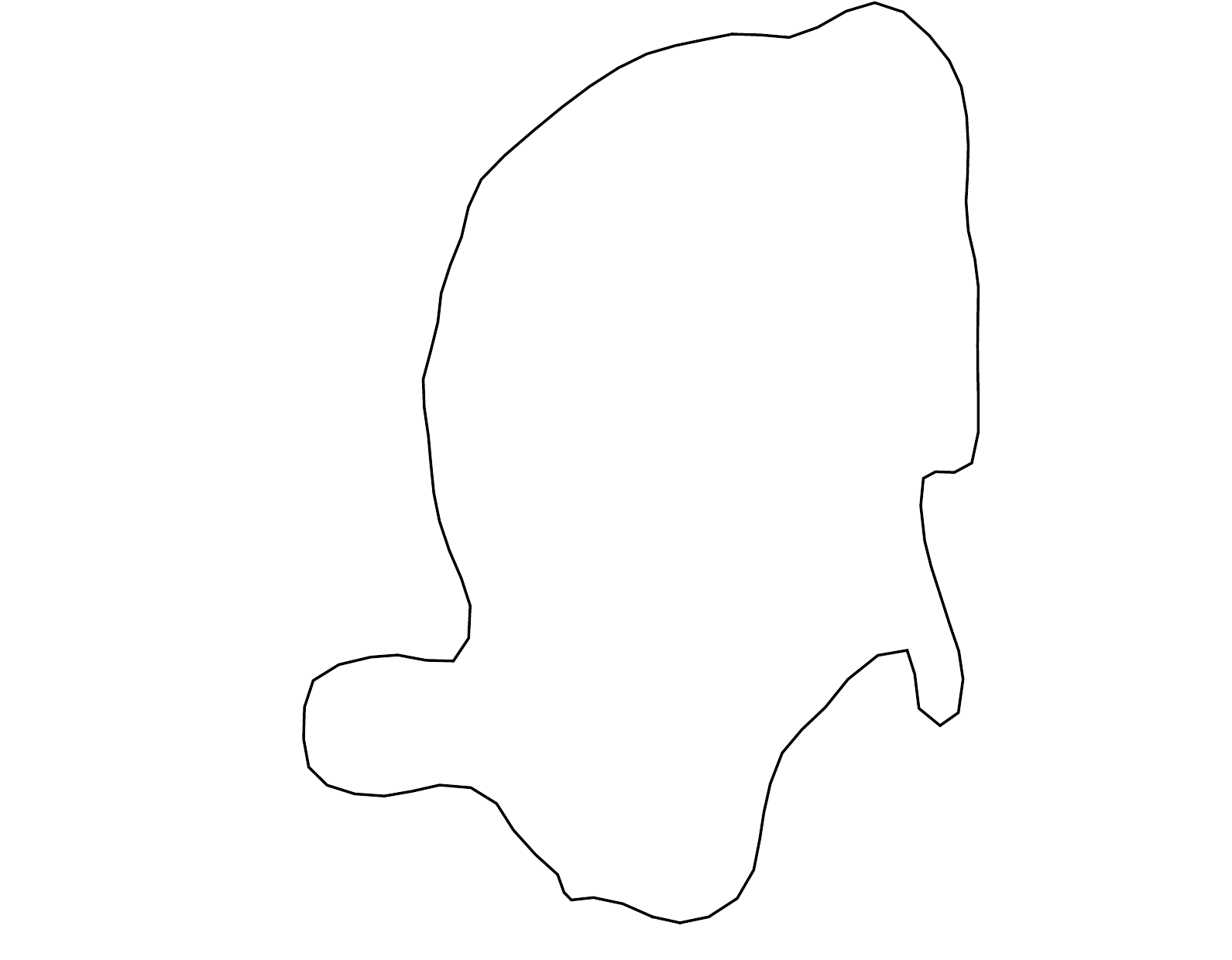}
\includegraphics[width=0.1\textwidth,angle=270]{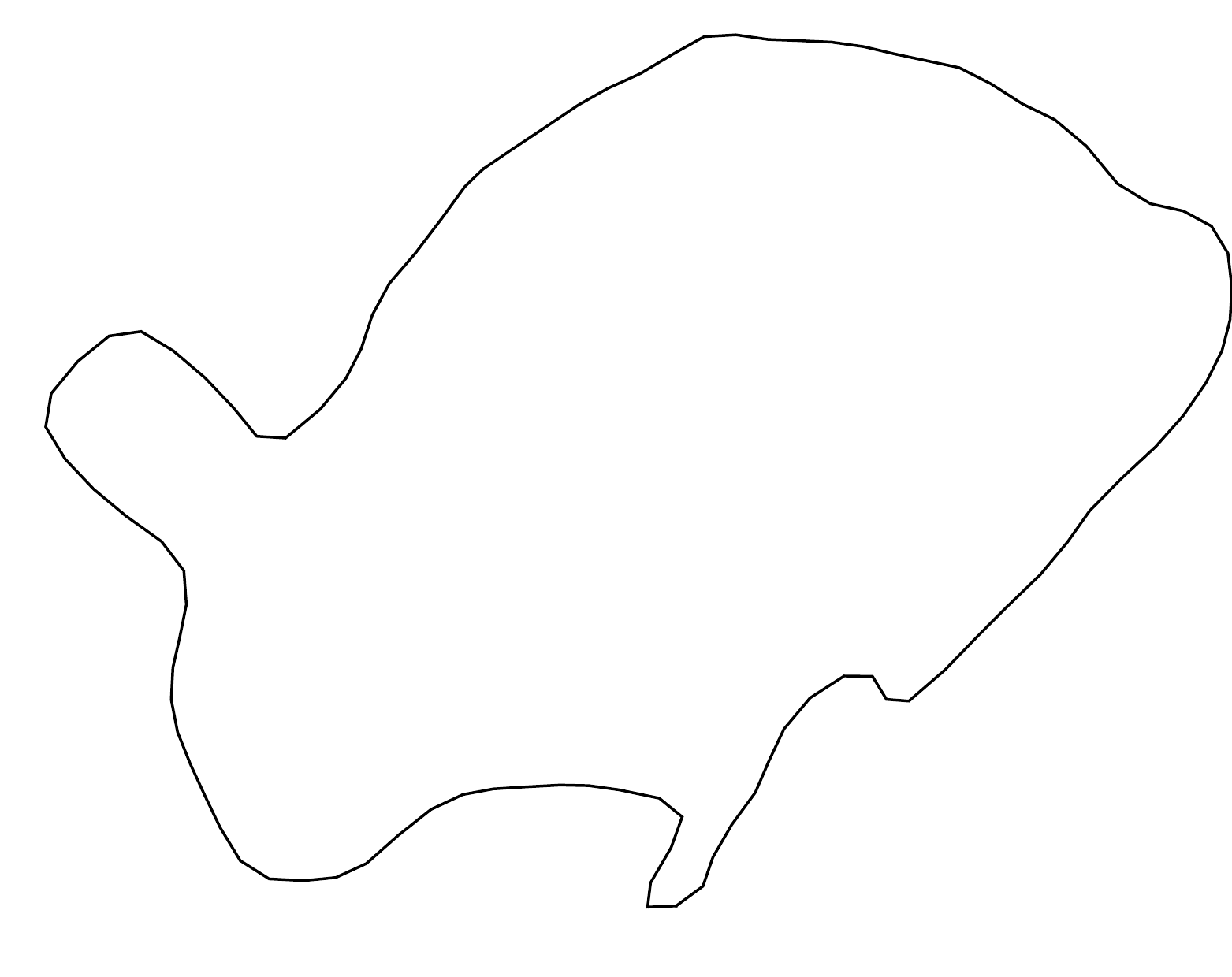}
\includegraphics[width=0.1\textwidth,angle=270]{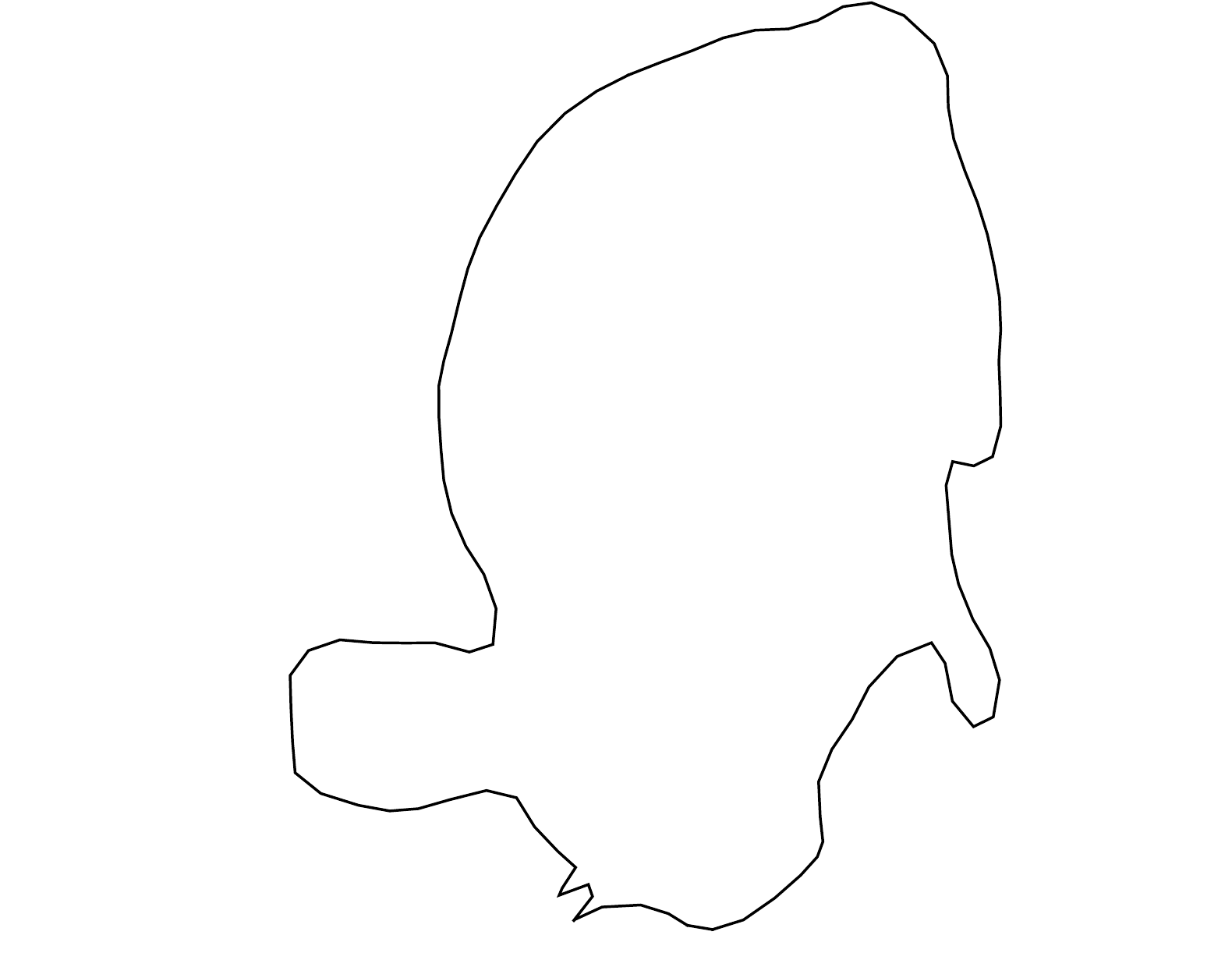}
\includegraphics[width=0.1\textwidth,angle=270]{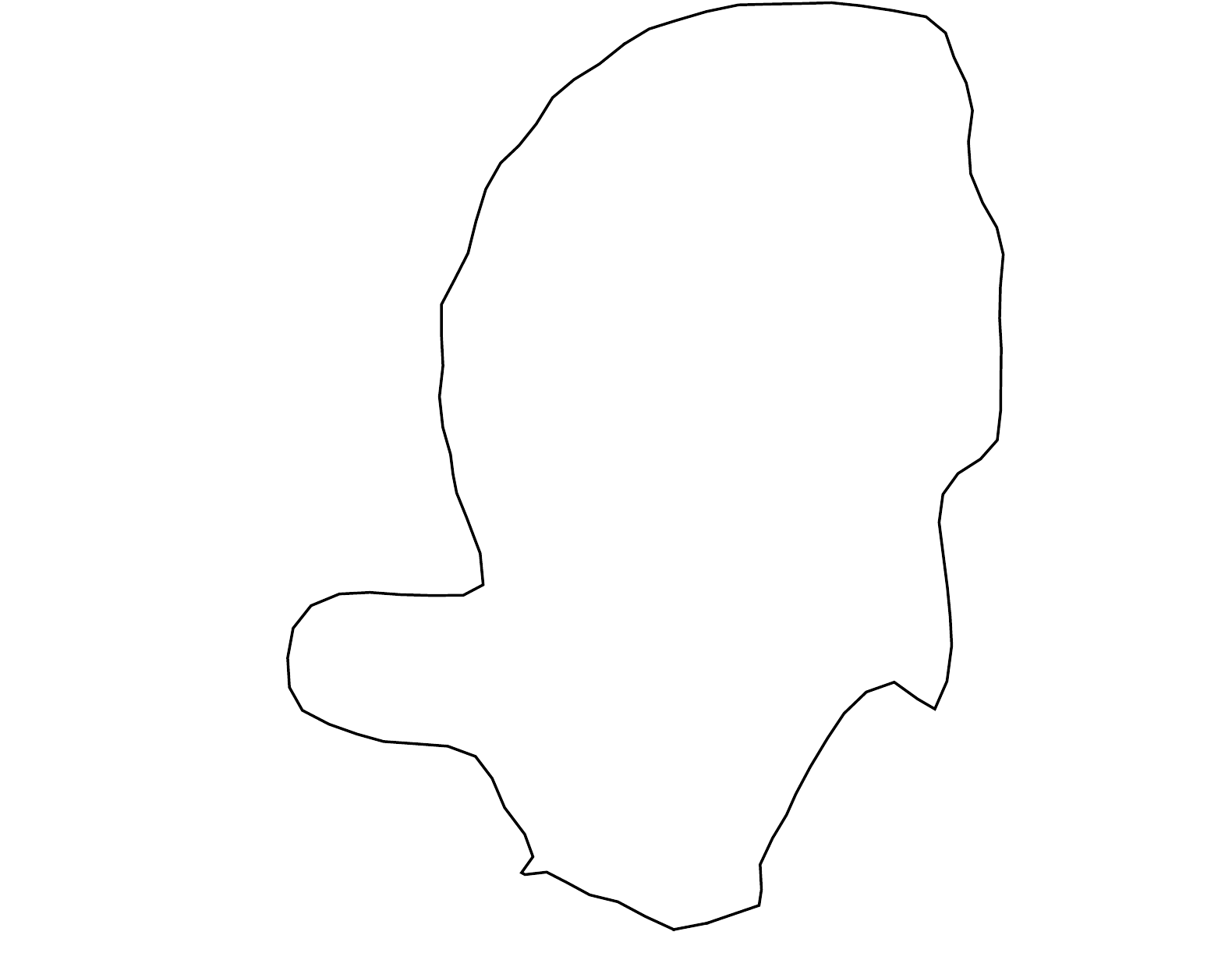}
\includegraphics[width=0.1\textwidth,angle=270]{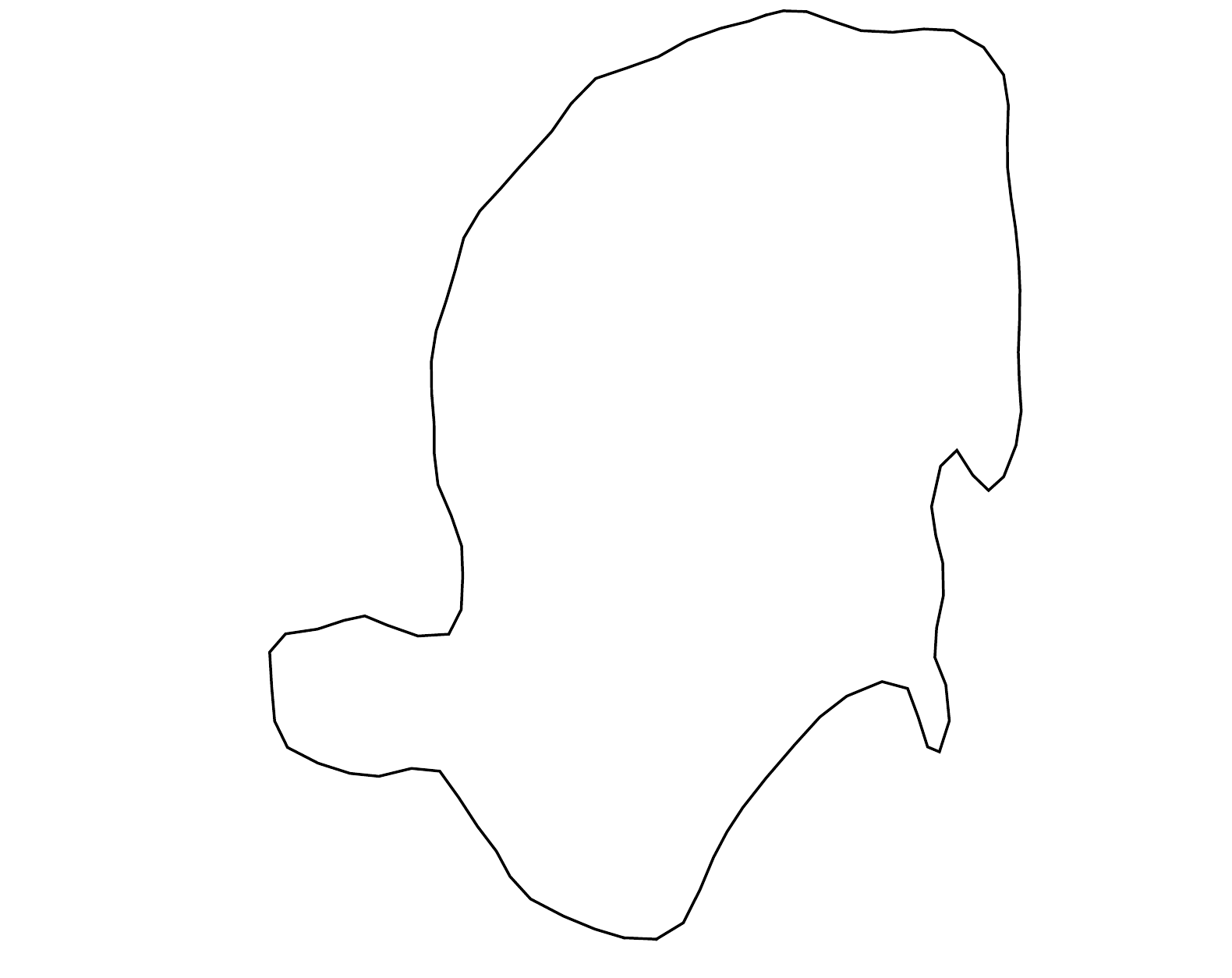}
\includegraphics[width=0.1\textwidth,angle=270]{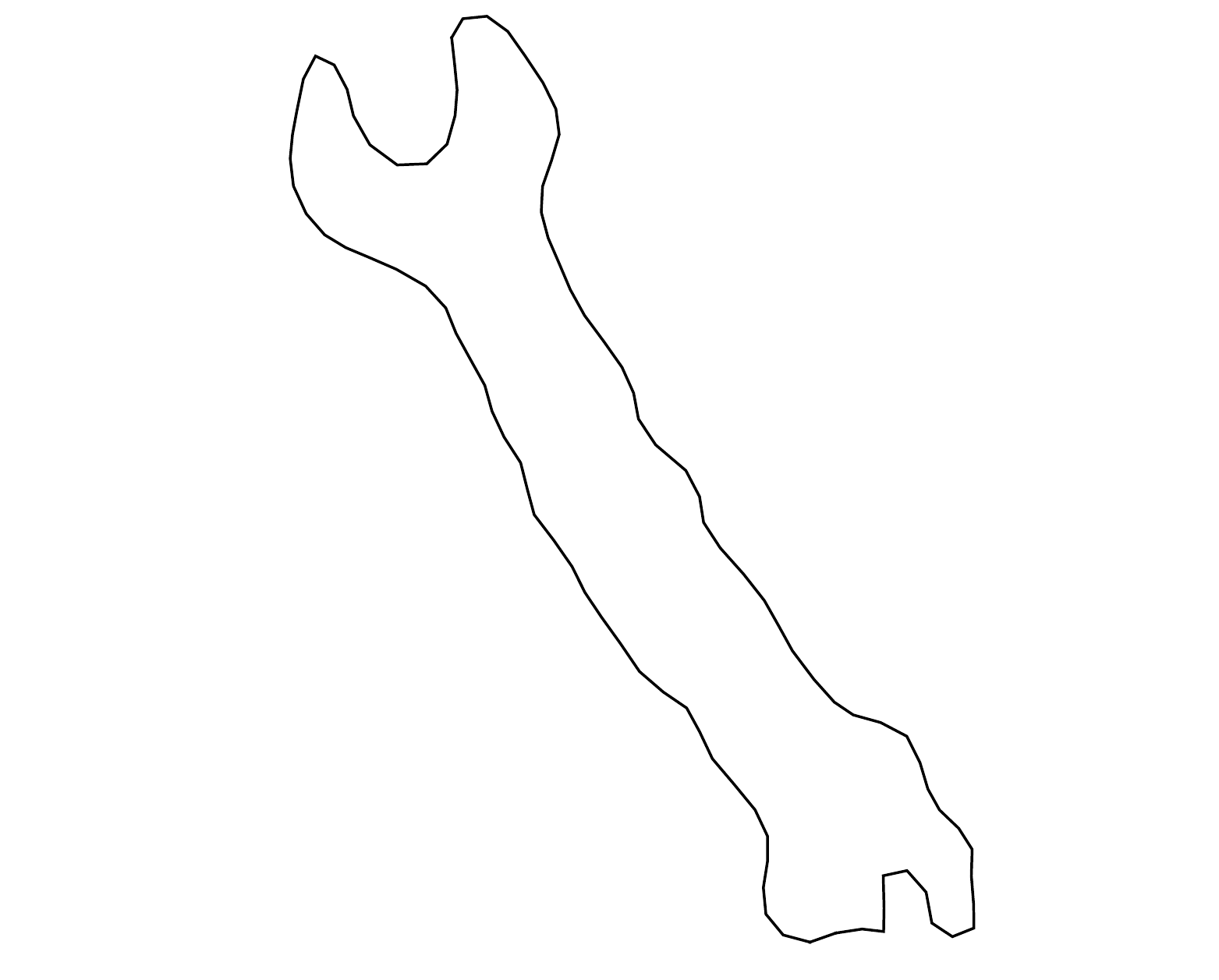}
\includegraphics[width=0.1\textwidth,angle=270]{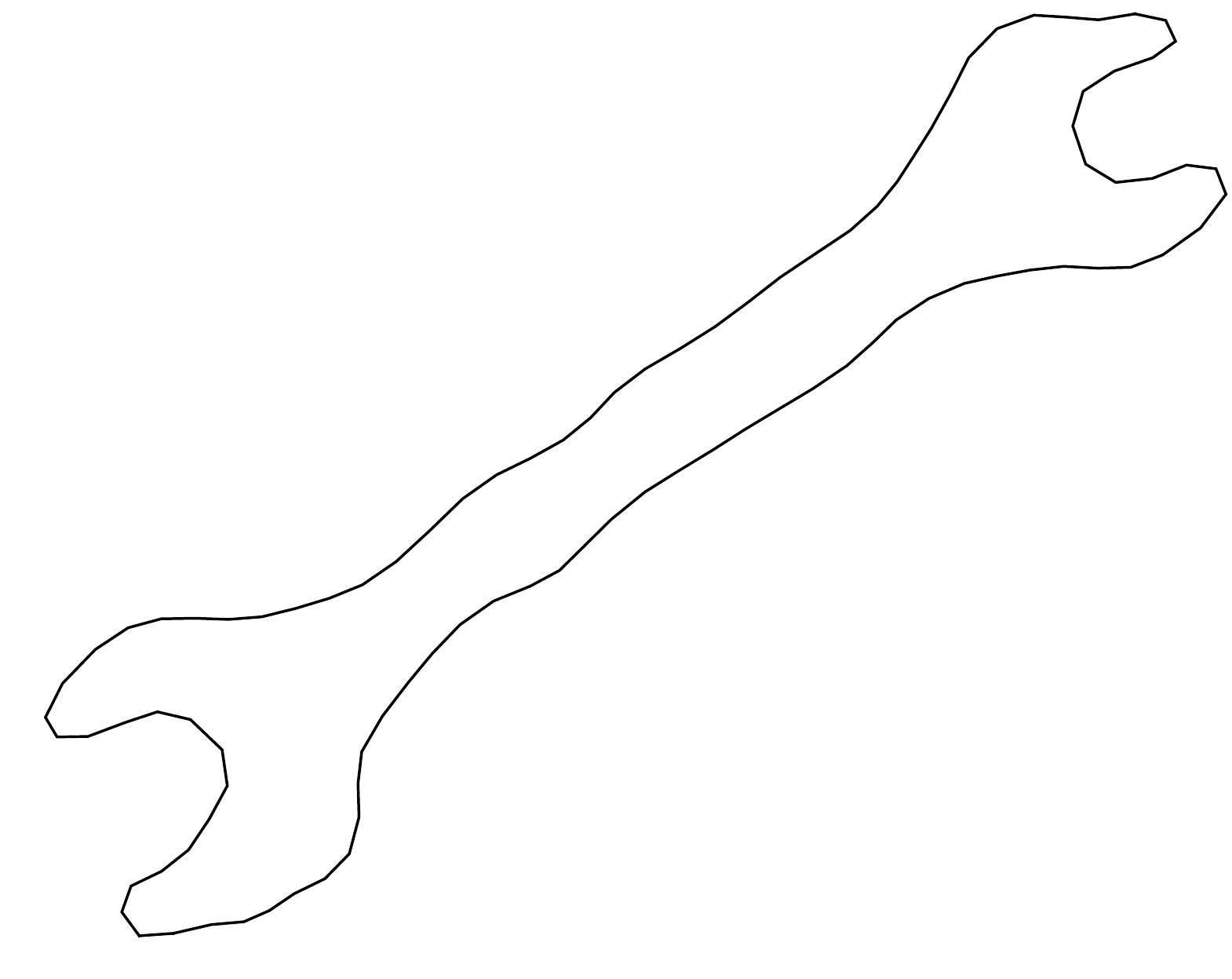}
\includegraphics[width=0.1\textwidth,angle=270]{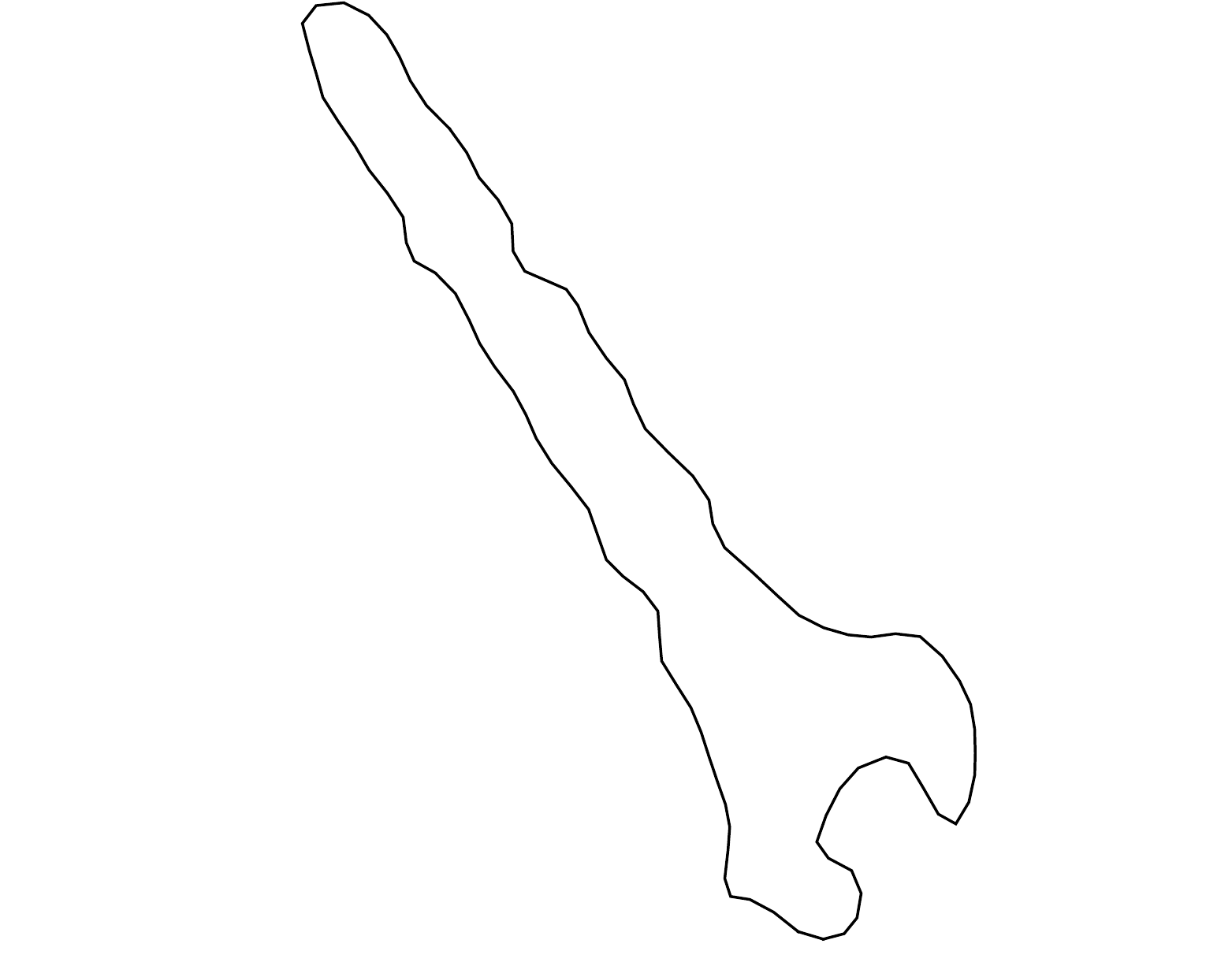}
\includegraphics[width=0.1\textwidth,angle=270]{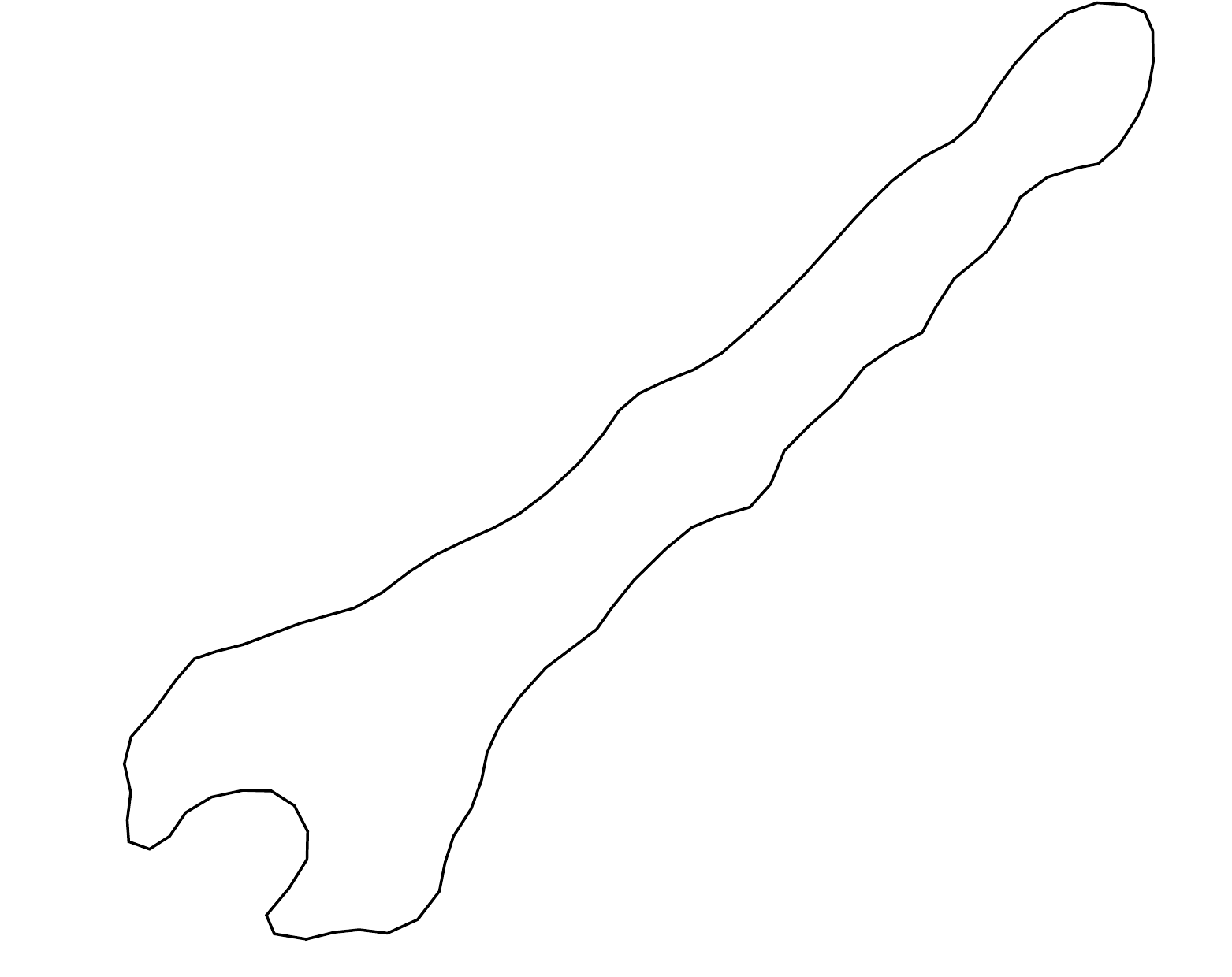}
\includegraphics[width=0.1\textwidth,angle=270]{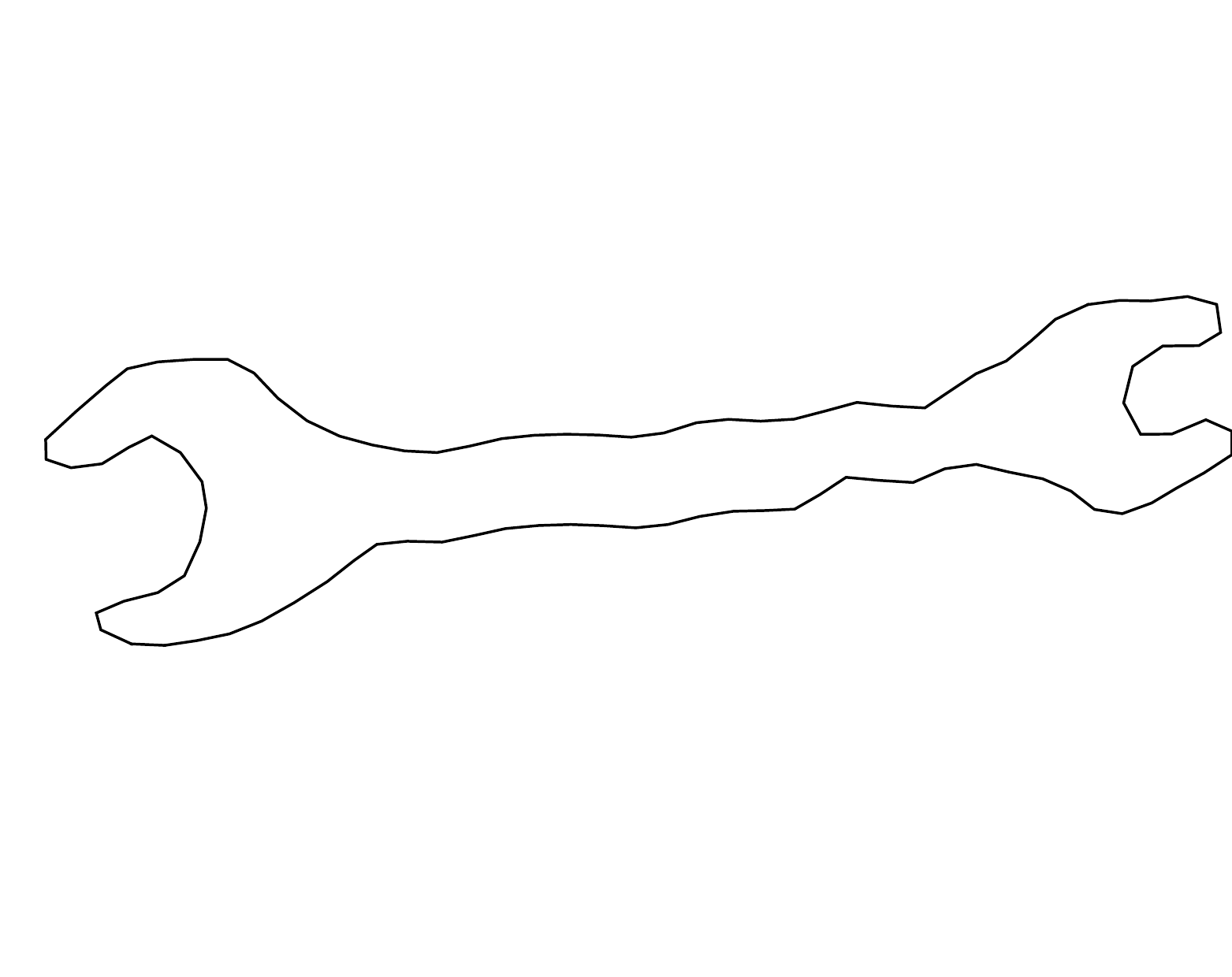}
\includegraphics[width=0.1\textwidth,angle=270]{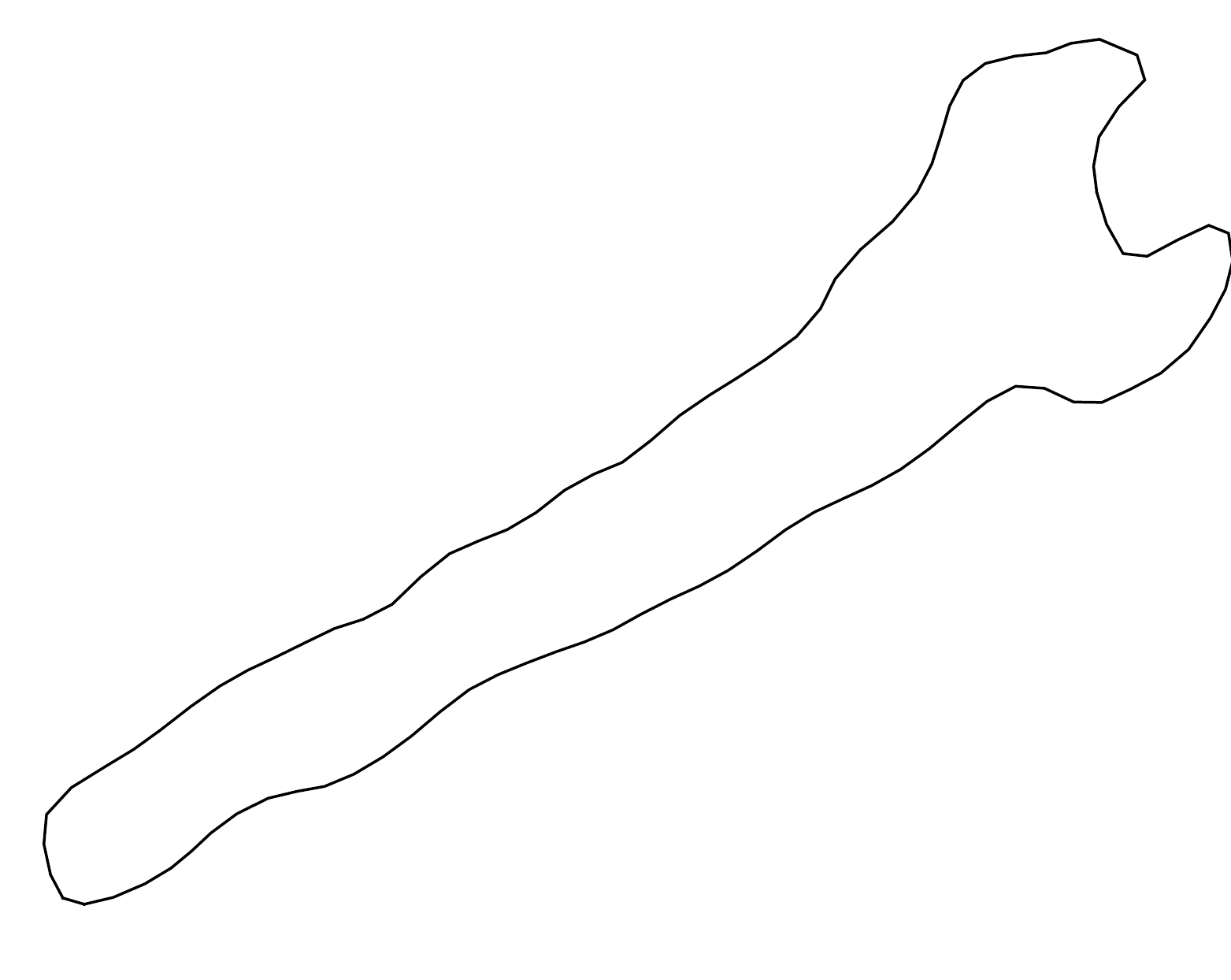}
\includegraphics[width=0.1\textwidth,angle=270]{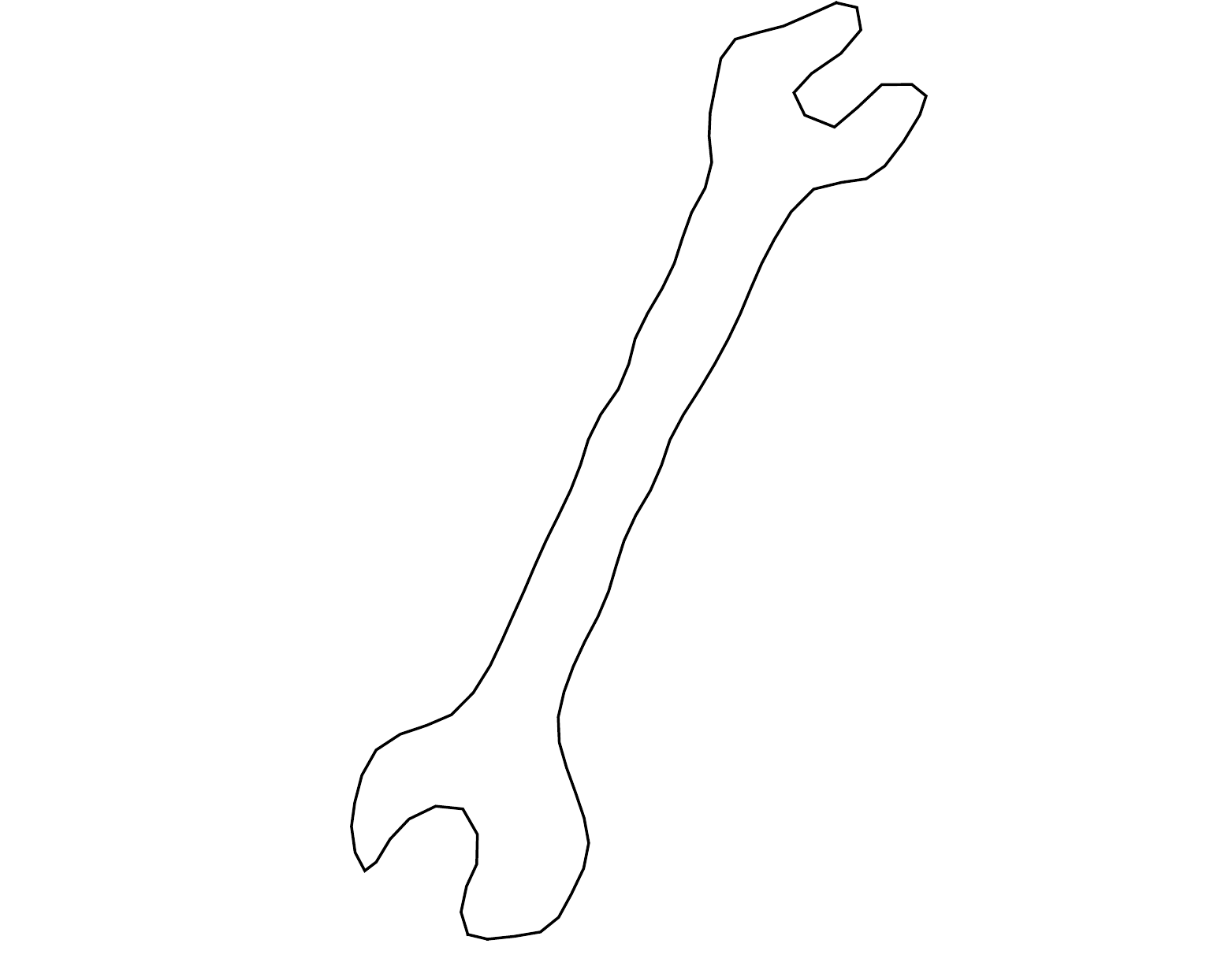}
\includegraphics[width=0.1\textwidth,angle=180]{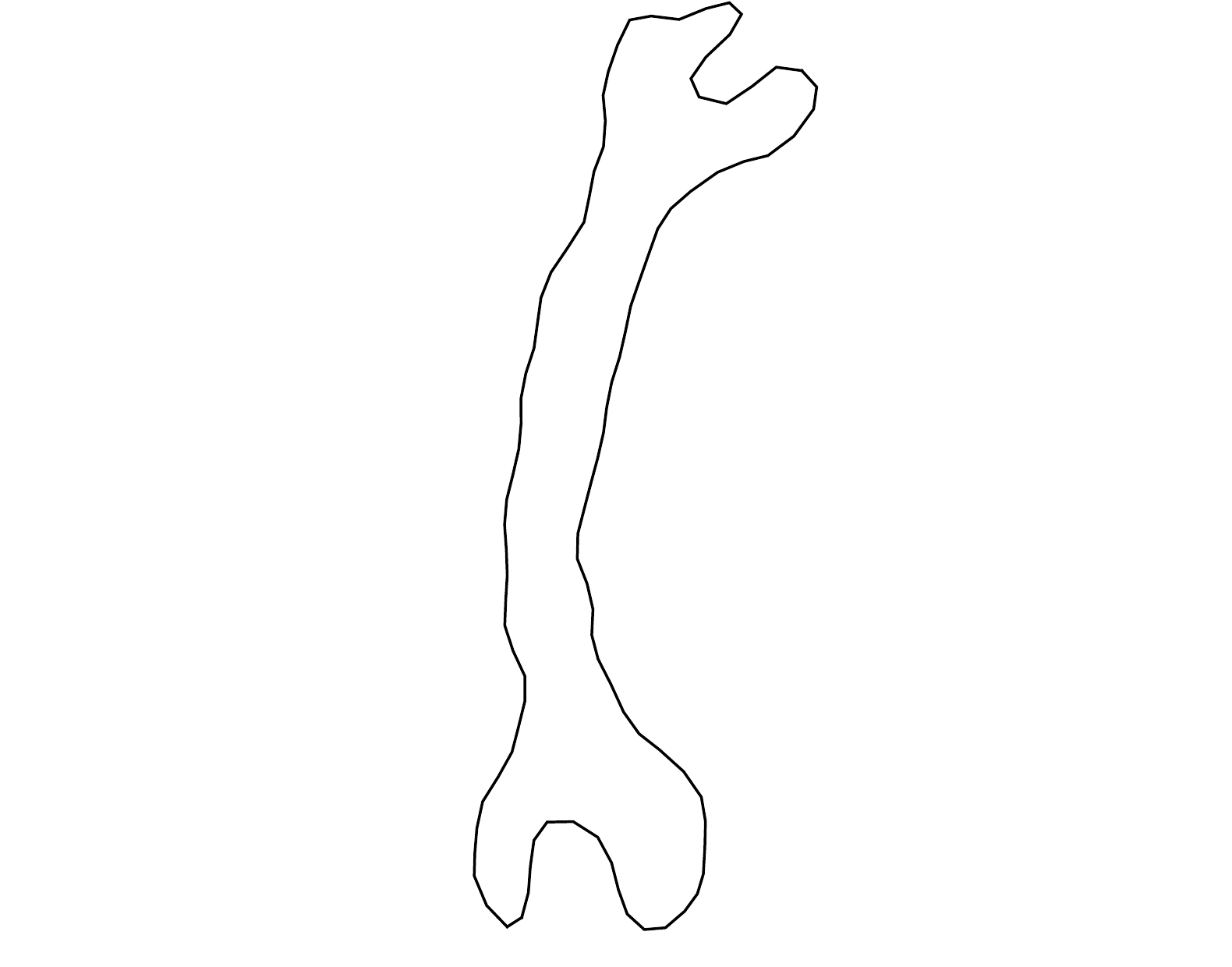}
\caption{Selection of shapes from the dataset \cite{kimia}.}
\label{fig:shapes}
\end{figure}

We tested our implementation on a dataset of shapes collected by the Computer Vision Group at Brown university \cite{kimia}. The dataset consists of black and white images of physical objects. It is natural to represent the objects by their boundaries using unparametrized curves. In addition to factoring out reparametrizations, we also factor out translations and rotations because we are not interested in the position of the curves in space. Some of the resulting curves are depicted in Fig.~\ref{fig:shapes}. We used splines of degree $n_\th=n_\phi=3$ with $N_\th=60$ and $N_\phi=20$ controls in space and of degree $n_t=2$ with $N_t=20$ controls in time. 

To solve the boundary value for geodesics as described in Sect. \ref{sec:boundary}, one has to construct an initial homotopy between the given boundary curves. If the curves differ only by a small deformation, this is unproblematic because the linear path connecting the curves can be used. For larger deformations, however, the linear path often has self-intersections, changing the winding number, which can lead to convergence problems during the energy minimization. Our solution was to construct homotopies by deforming the initial curve to a circle and the circle to the target curve using linear paths. This worked well for all examples presented in this article. It also led to the same optima in all cases where the linear path was free of self-intersections. In future work, we plan to investigate further ways of creating good and natural initial paths. 

\begin{figure}
\centering
	\includegraphics[width=.2\textwidth]{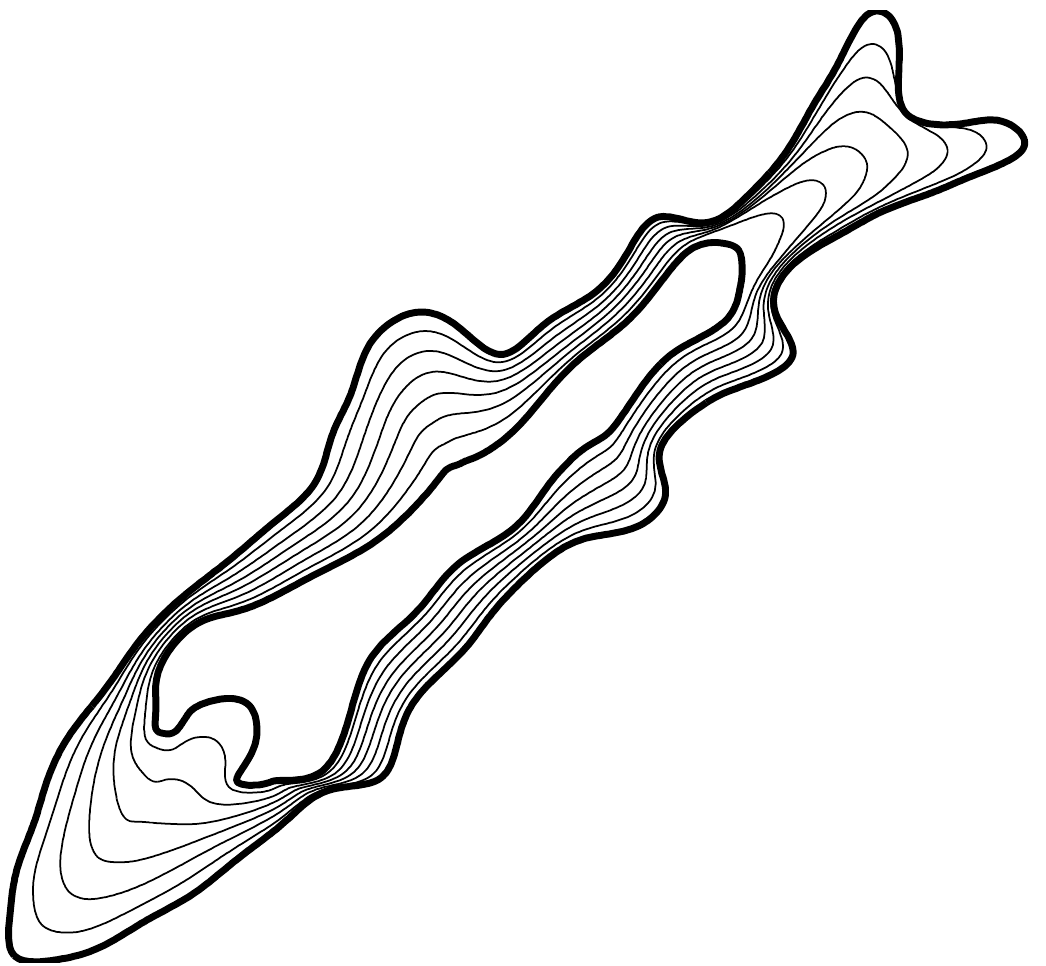}
	\includegraphics[width=.2\textwidth]{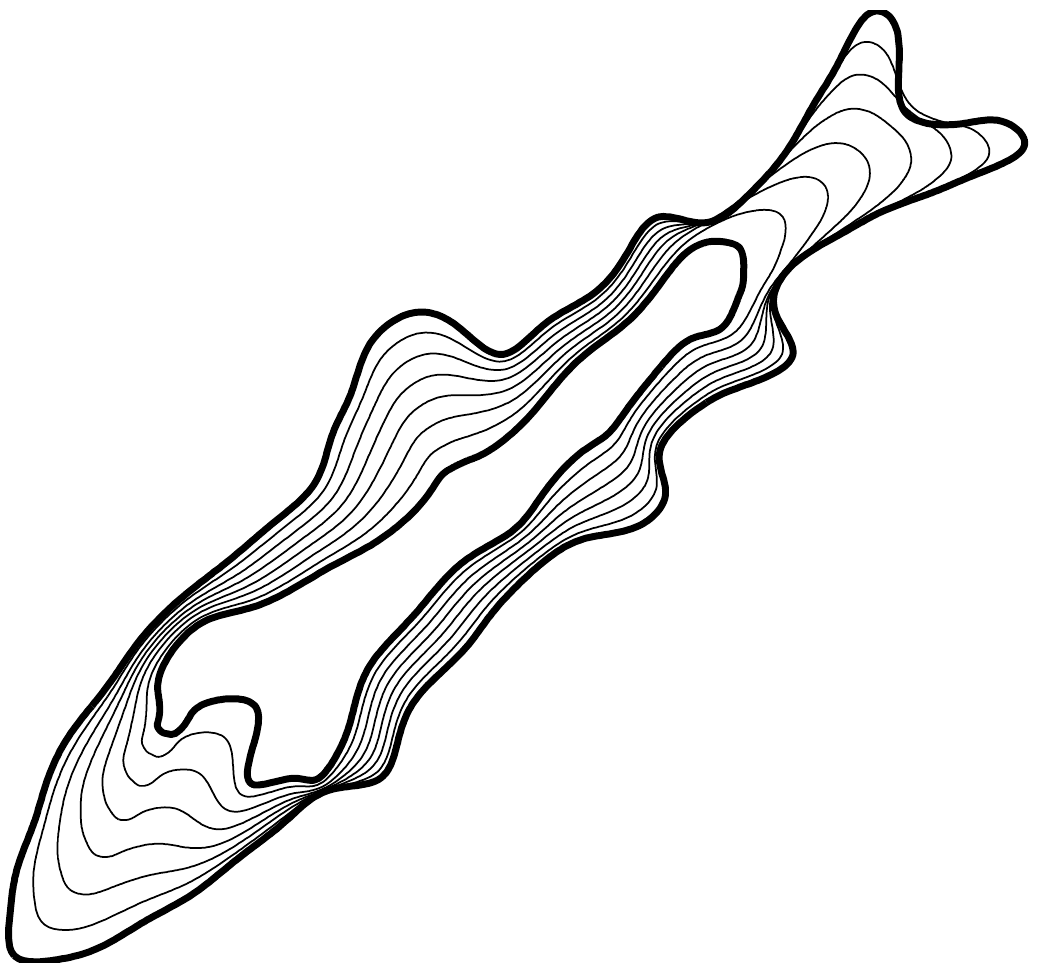}	
	\includegraphics[width=.2\textwidth]{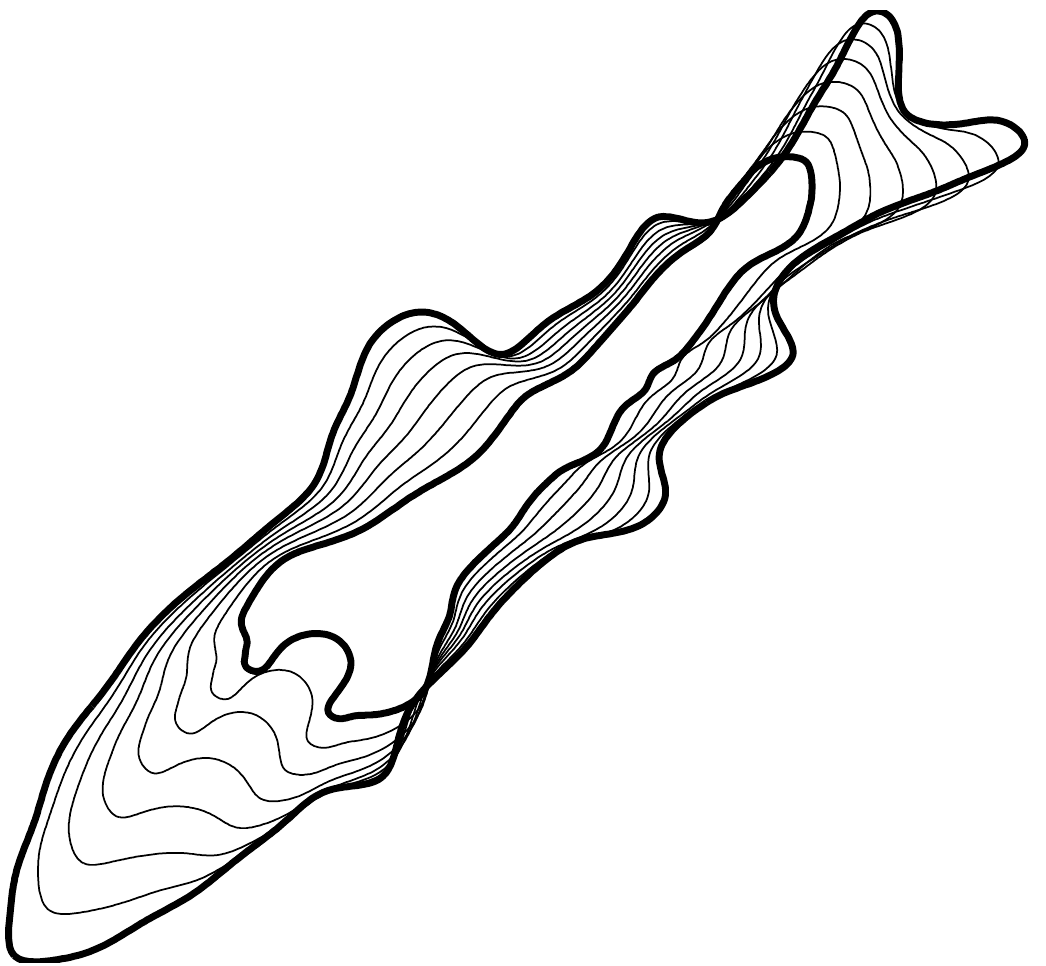}		
	\includegraphics[width=.2\textwidth]{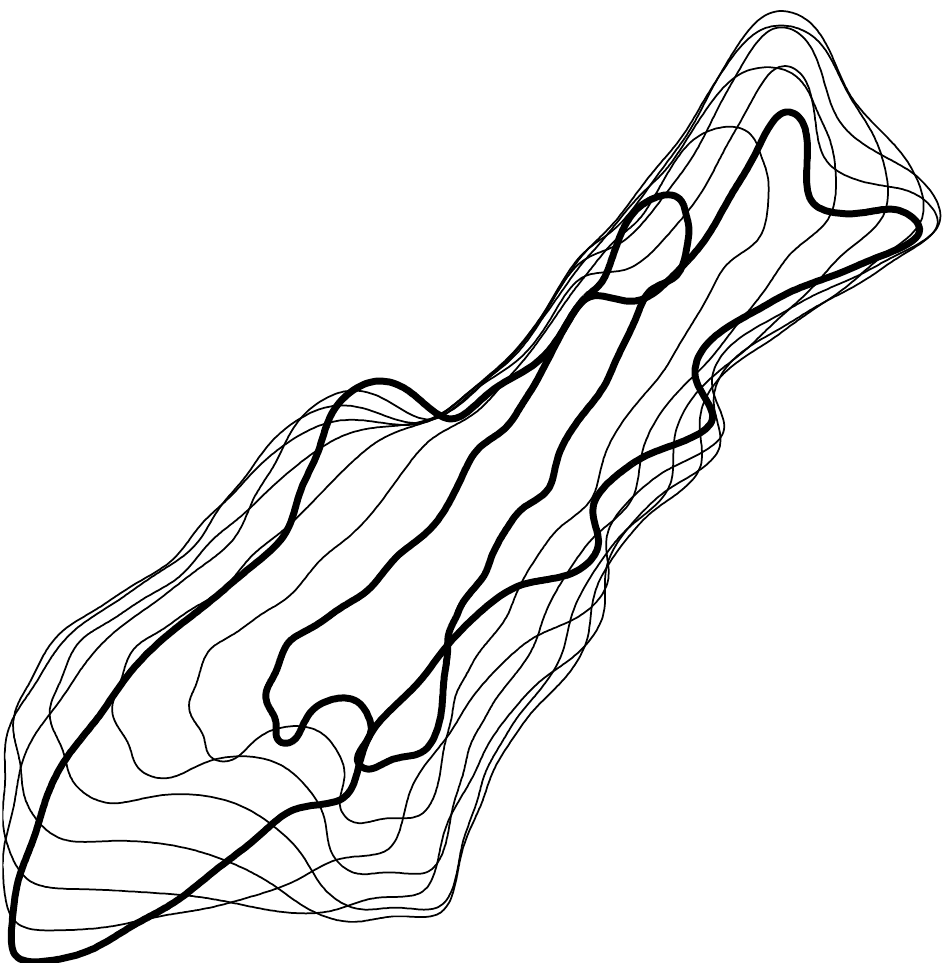}
	\caption{Geodesics between a fish and a tool in the space of unparametrized curves. The metric parameter $a_2$ is increased by a factor 10 in the second, a factor 100 in the third, and a factor 1000 in the fourth column. The corresponding geodesic distances are $138.02$, $162.55$, $246.46$ and $468.74$. Note that since we also optimize over translations and rotations of the target curve, the position in space varies.}
	\label{fig:geodesicsConstants}
\end{figure}

The choice of parameters $a_0$, $a_1$, and $a_2$ of the Riemannian metric can have a large influence on the resulting optimal deformations, as can be seen in Fig.~\ref{fig:geodesicsConstants}. In this article, we used the following ad-hoc strategy for choosing the constants: we computed the average $L^2$-, $H^1$- and $H^2$-contributions $\overline{E}_0$, $\overline{E}_1$, $\overline{E}_2$ to the energy of linear paths between each pair of curves in the dataset. Then we normalized $a_0$ to $1$ and chose $a_1$ and $a_2$ such that
\[
a_0 \overline{E}_0 : a_1 \overline{E}_1 : a_2 \overline{E}_2 = 1 : 1 : 1\;\text{ and }\;
\overline{E} = a_0 \overline{E}_0 + a_1 \overline{E}_1 +
a_2 \overline{E}_2 = 100\,.
\]
This resulted in the parameter values $a_1 = 250$ and $a_2 = 0.004$. It remains open how to best choose the constants depending on the data under consideration.

\subsection{Comparison to non-elastic metrics}

Riemannian metrics on spaces of unparametrized curves are often called elastic metrics, as they allow both for bending and stretching of the curve. In the elastic case, solving the boundary value problem for geodesics involves optimizing over the $\operatorname{Diff}(S^1)$-orbit of the initial or final shape. This is a computationally expensive and difficult task since the  diffeomorphism group is infinite-dimensional. 

An alternative and simpler approach is to parametrize the curves by unit-speed and to calculate geodesics in the space of parametrized curves. This could be called a non-elastic approach. (Of course, one might wish to also factor out rigid transformations and constant shifts of the parametrization, but this is much simpler because these groups are finite dimensional.) 

\begin{figure}
\centering
	\includegraphics[width=.2\textwidth,angle=270]{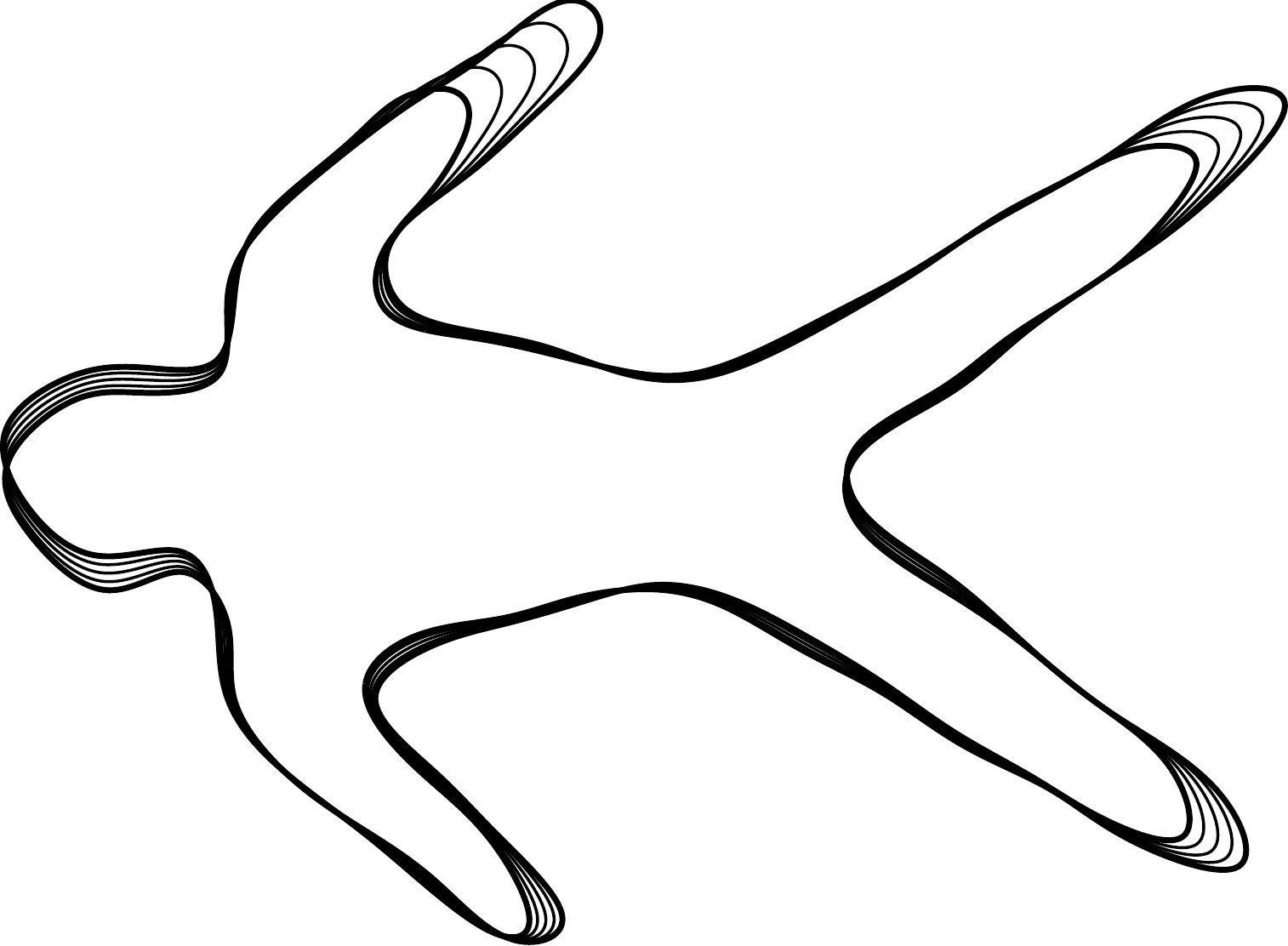}
		\hspace{1cm}
	\includegraphics[width=.2\textwidth,angle=270]{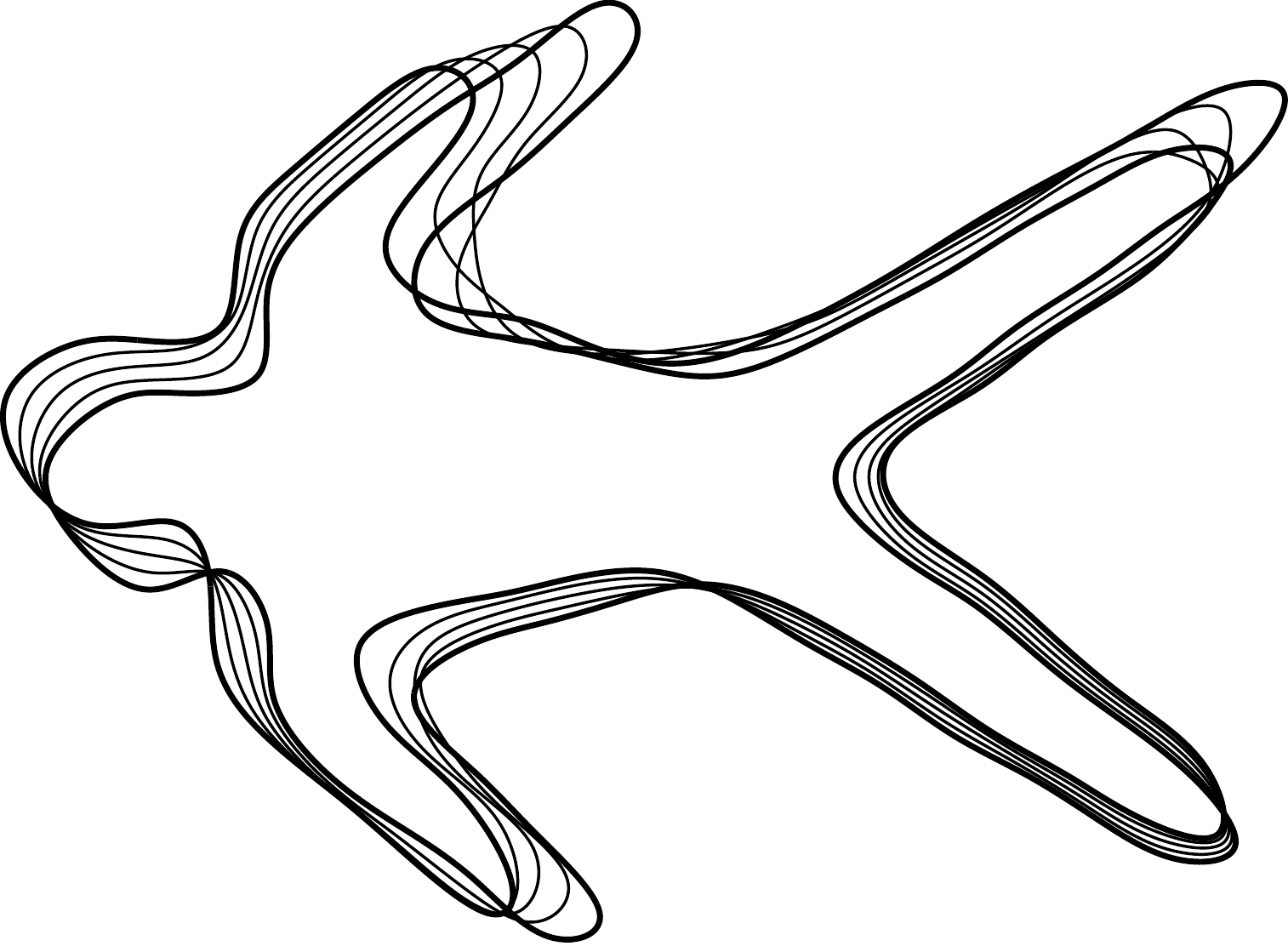}
		\hspace{1cm}
	\includegraphics[width=.2\textwidth,angle=270]{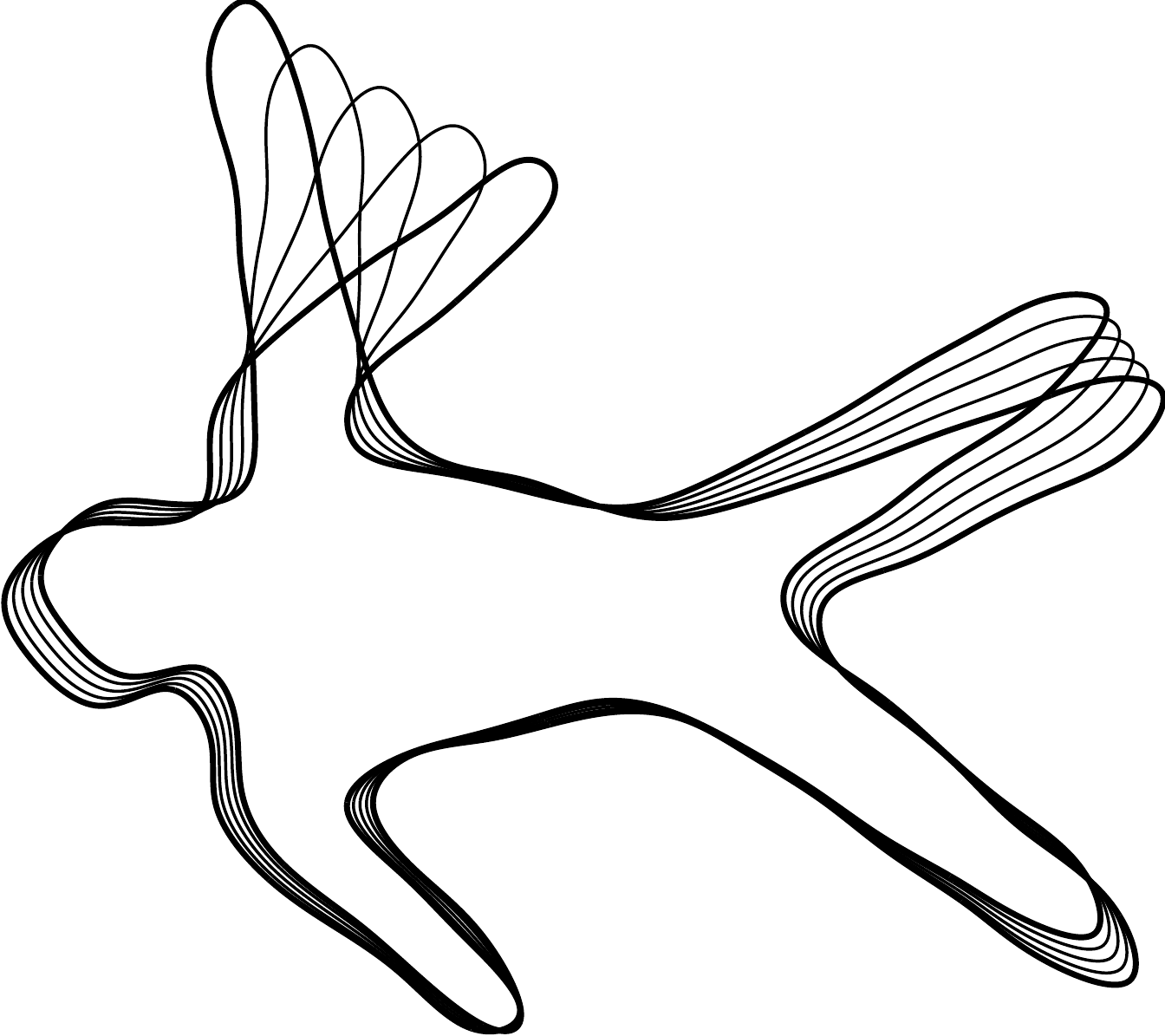}
		\hspace{1cm}
	\includegraphics[width=.2\textwidth,angle=270]{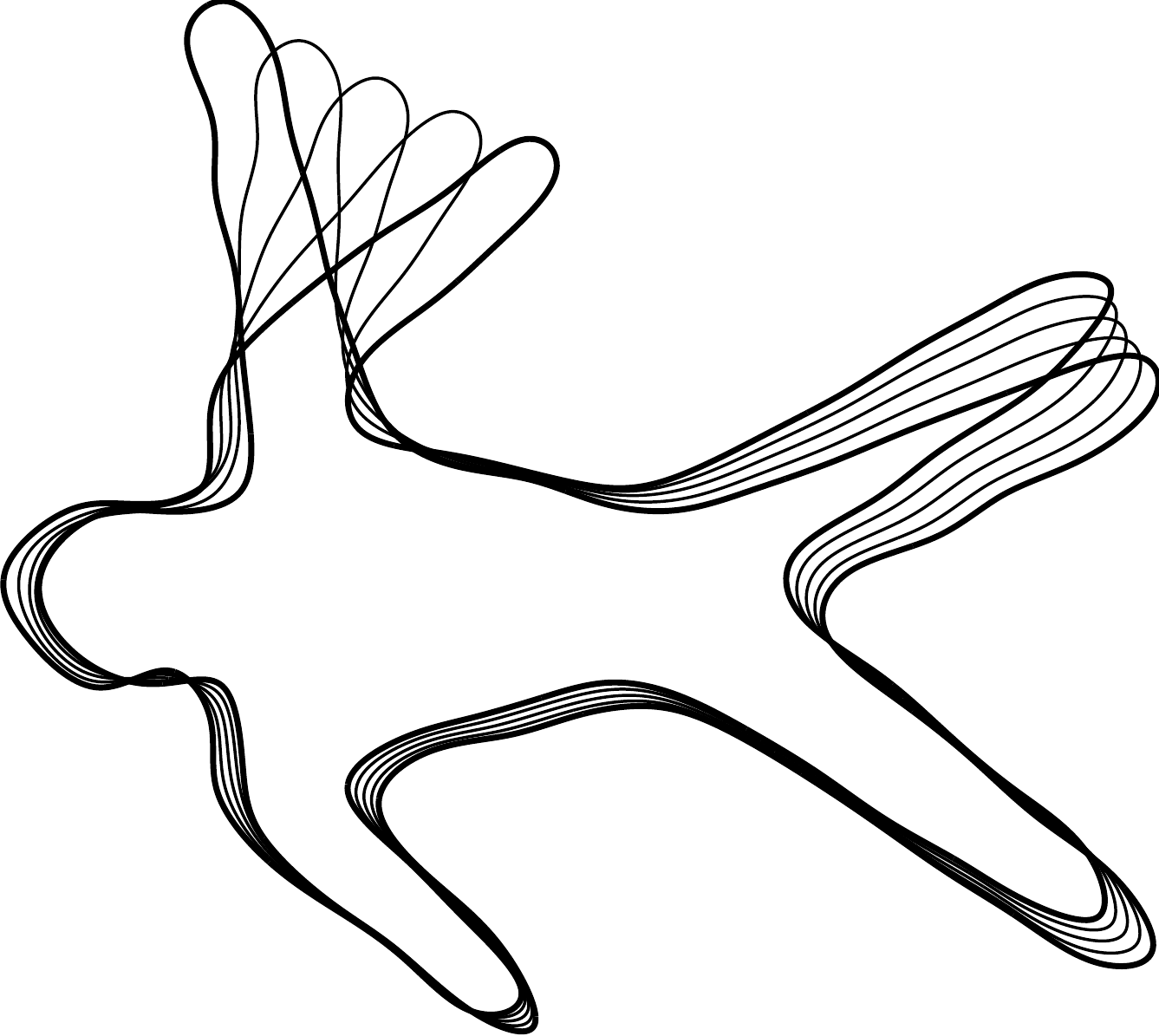}
	\caption{Geodesics in the space of unparametrized (1st, 3rd) versus parametrized (2nd, 4th) curves modulo rotations and translations. Note that since we also optimize over translations and rotations of the target curve, the curves in the first two panels are  aligned differently.}
	\label{fig:geodesics}
\end{figure}

Our experiments suggest that using the more involved elastic approach  pays off. The two approaches yield different results, and in particular in cases involving large amounts of stretching the geodesics found using the first approach appear more natural. However, in cases involving mainly bending of the curves the results are very similar (see Fig.~\ref{fig:geodesics}). 

\subsection{Clustering and principal component analysis}

\begin{figure}
\centering
	\includegraphics[width=.4\textwidth]{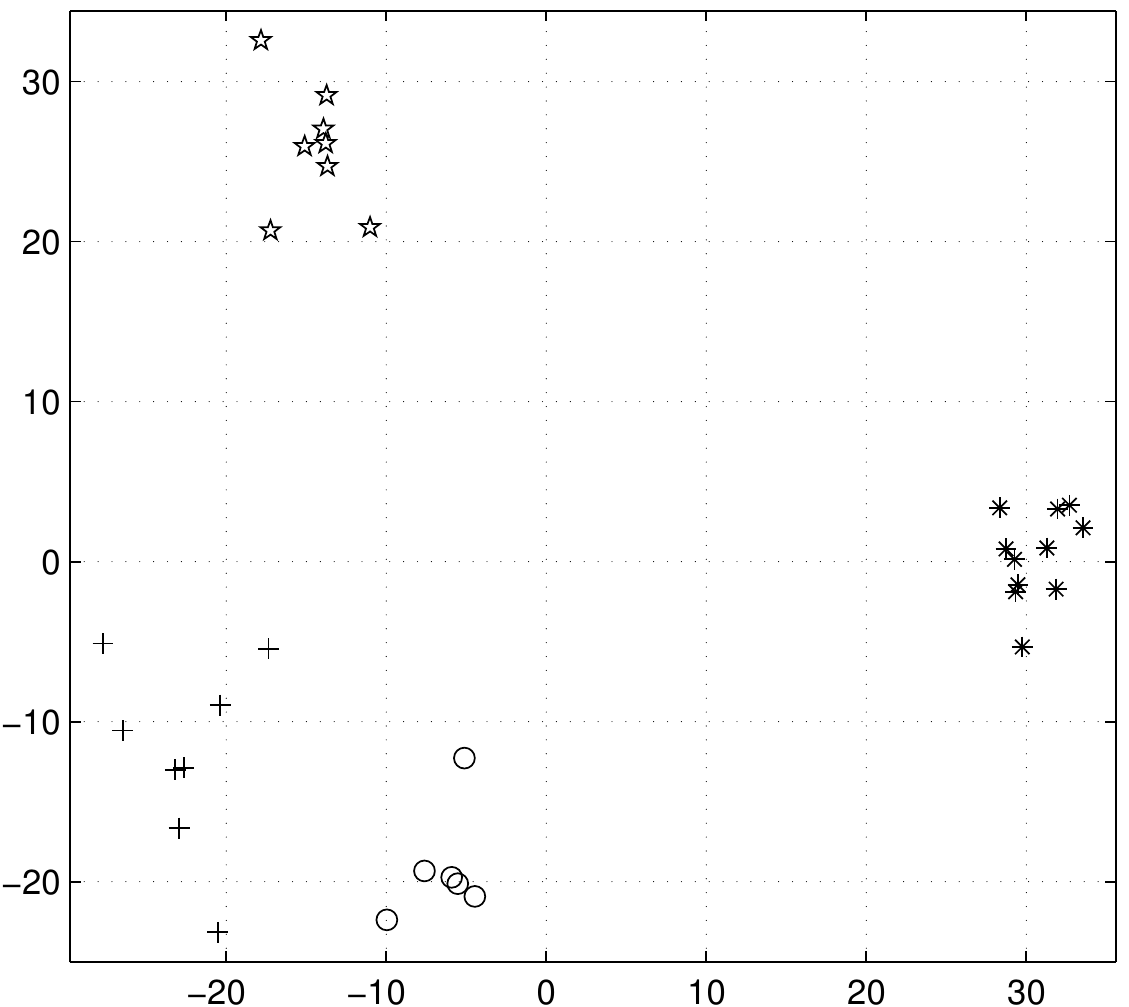}
		\hspace{1cm}
	\includegraphics[width=.4\textwidth]{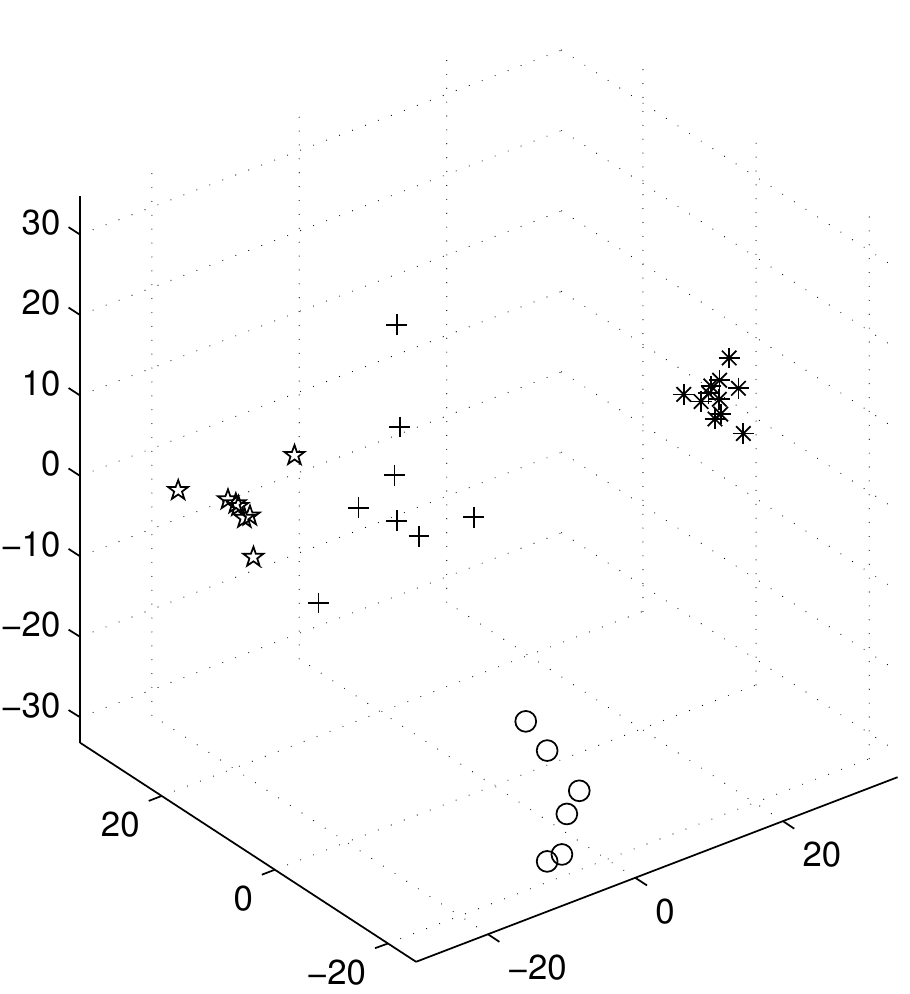}
	\caption{The matrix of geodesic distances between the shapes in Fig.~\ref{fig:shapes}, visualized using multi-dimensional scaling in two and three dimensions. The labels are: fish $\medcircle$, sting rays $\medstar$, bunnies $\varhexstar$, tools $+$.}
	\label{fig:mds}
\end{figure}

We investigated if pairwise geodesic distances can be used to cluster shapes into meaningful groups. To this aim, we calculated all pairwise distances between the 33 shapes presented in Fig.~\ref{fig:shapes}. Solving the corresponding 528 boundary value problems took about two hours on a 3GHz processor with four cores. The resulting distance matrix is visualized in Fig.~\ref{fig:mds} using multi-dimensional scaling; the plot suggests that objects of the same group lie close together. Indeed, agglomerative clustering with 4 clusters reproduces exactly the subgroups of the database.

\begin{figure}
\centering
	\includegraphics[width=.20\textwidth,angle=275]{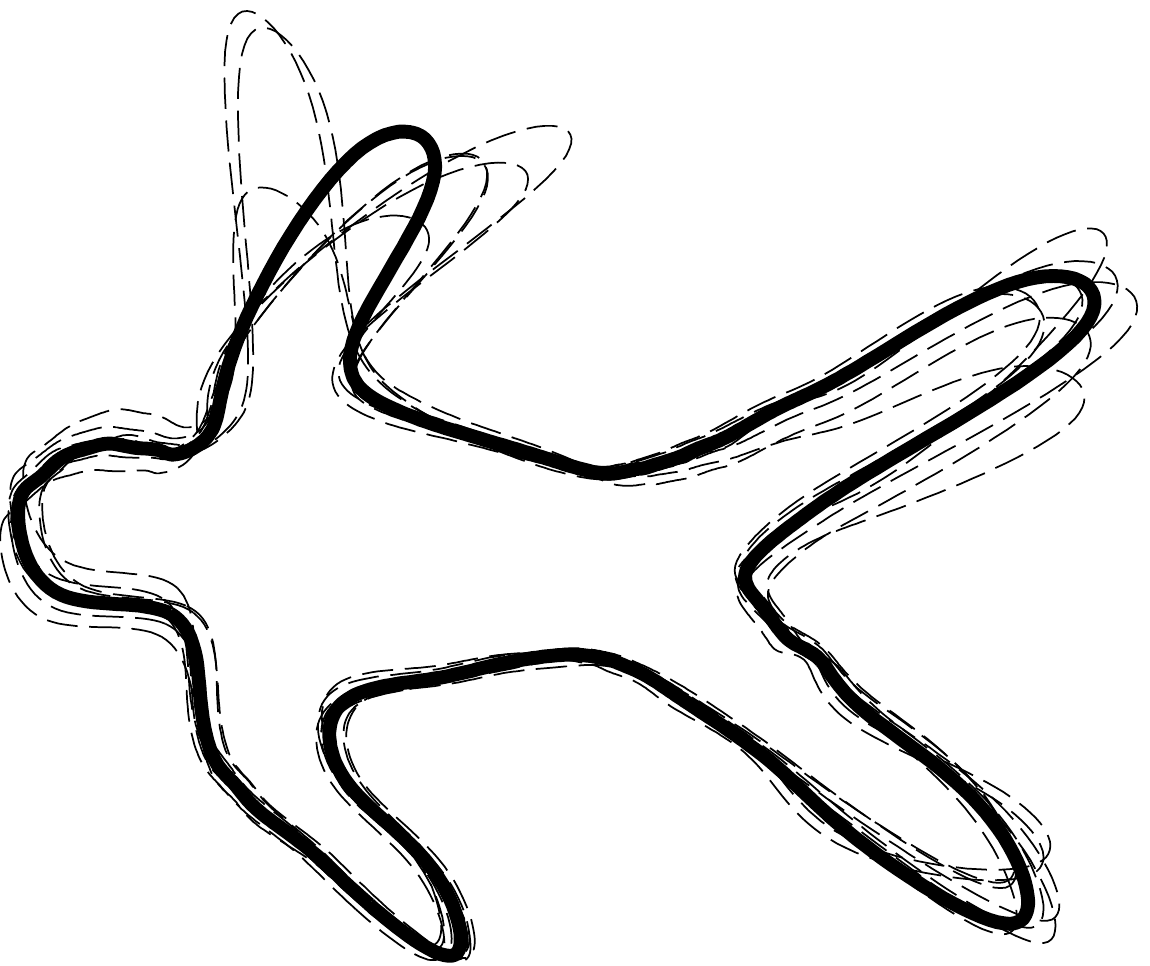}\hspace{1cm}
	\includegraphics[width=.20\textwidth,angle=275]{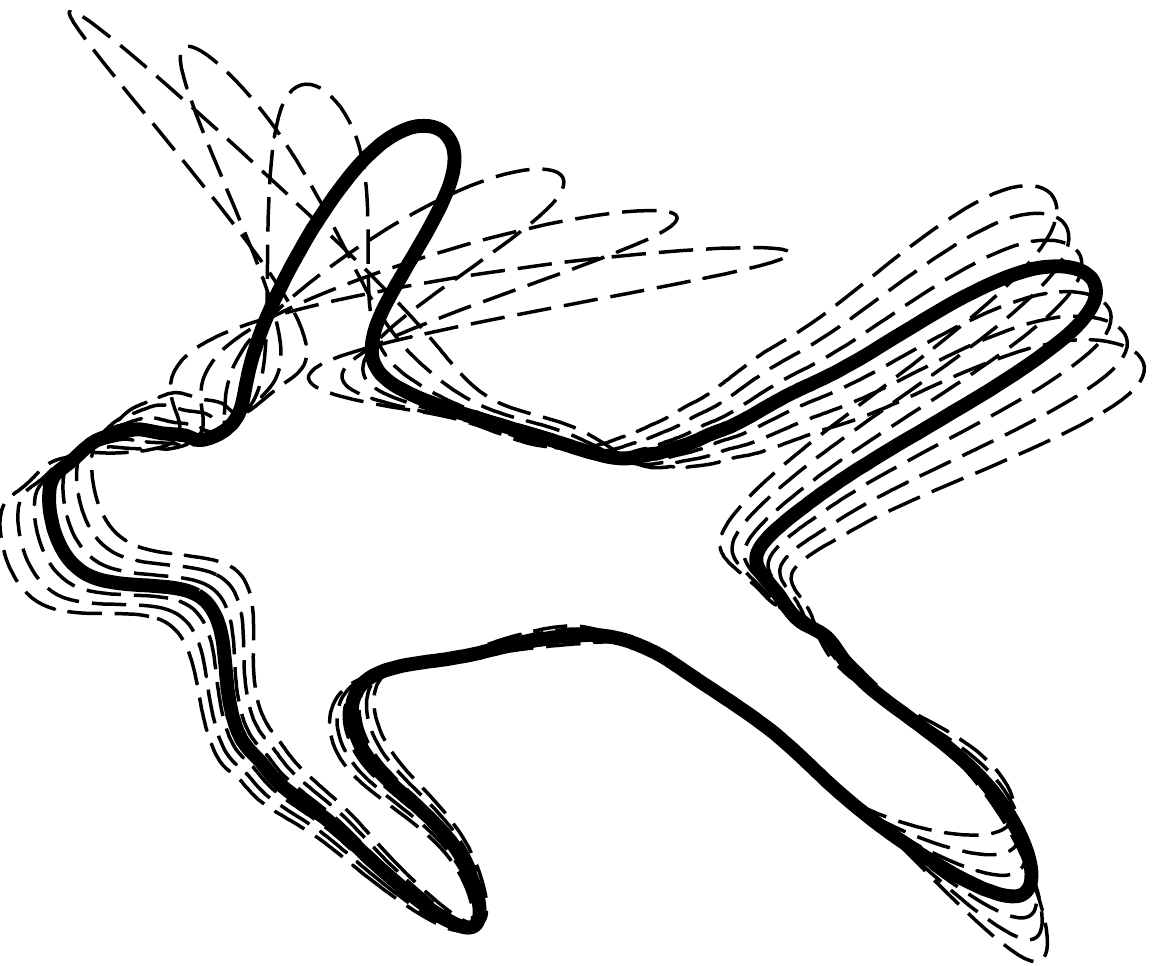}\hspace{1cm}
	\includegraphics[width=.20\textwidth,angle=275]{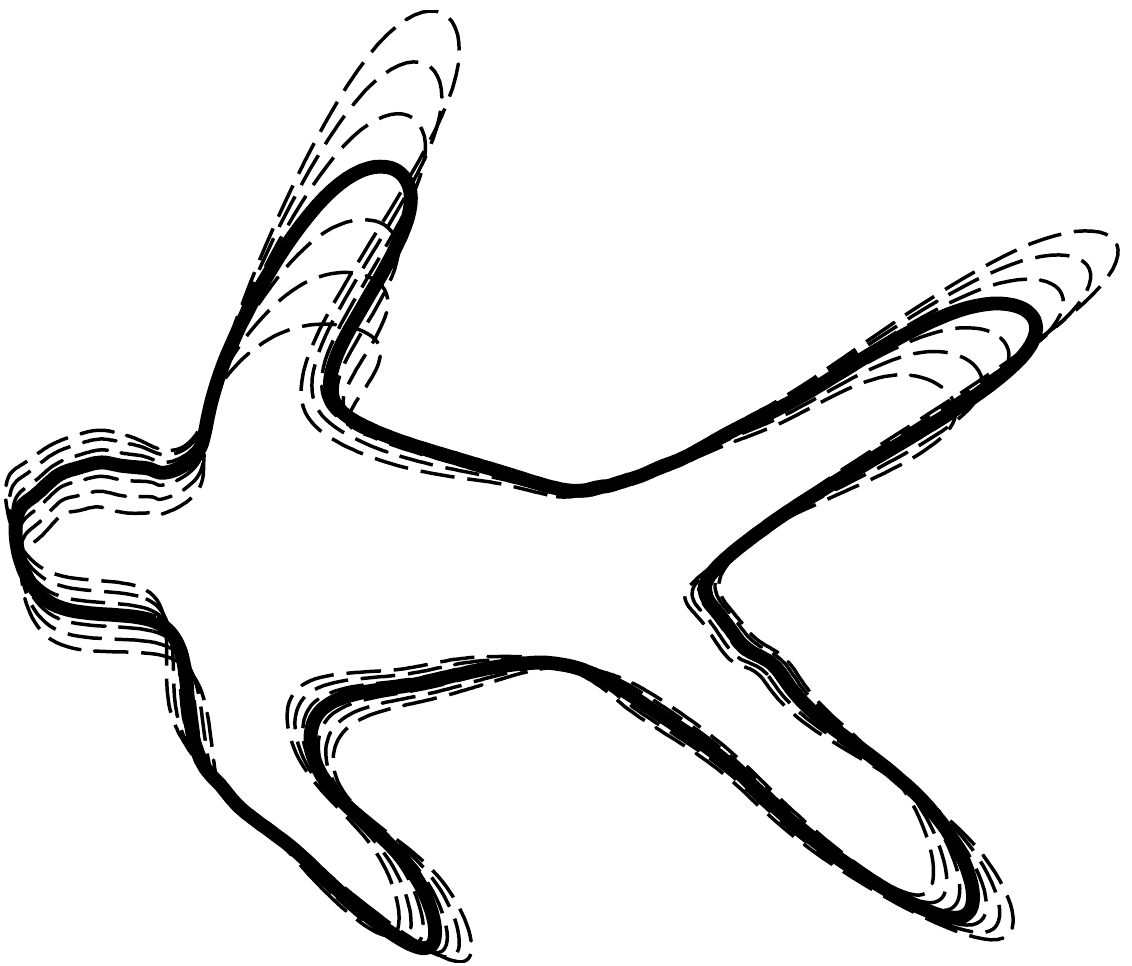}\\
	\includegraphics[width=.20\textwidth,angle=275]{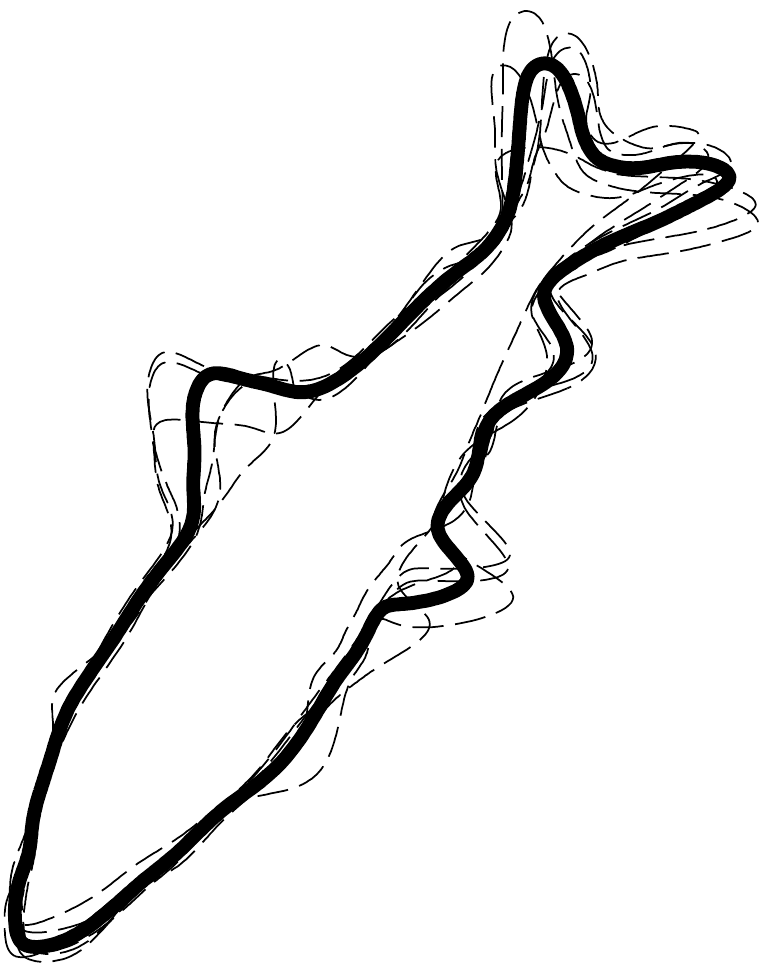}
	\includegraphics[width=.20\textwidth,angle=275]{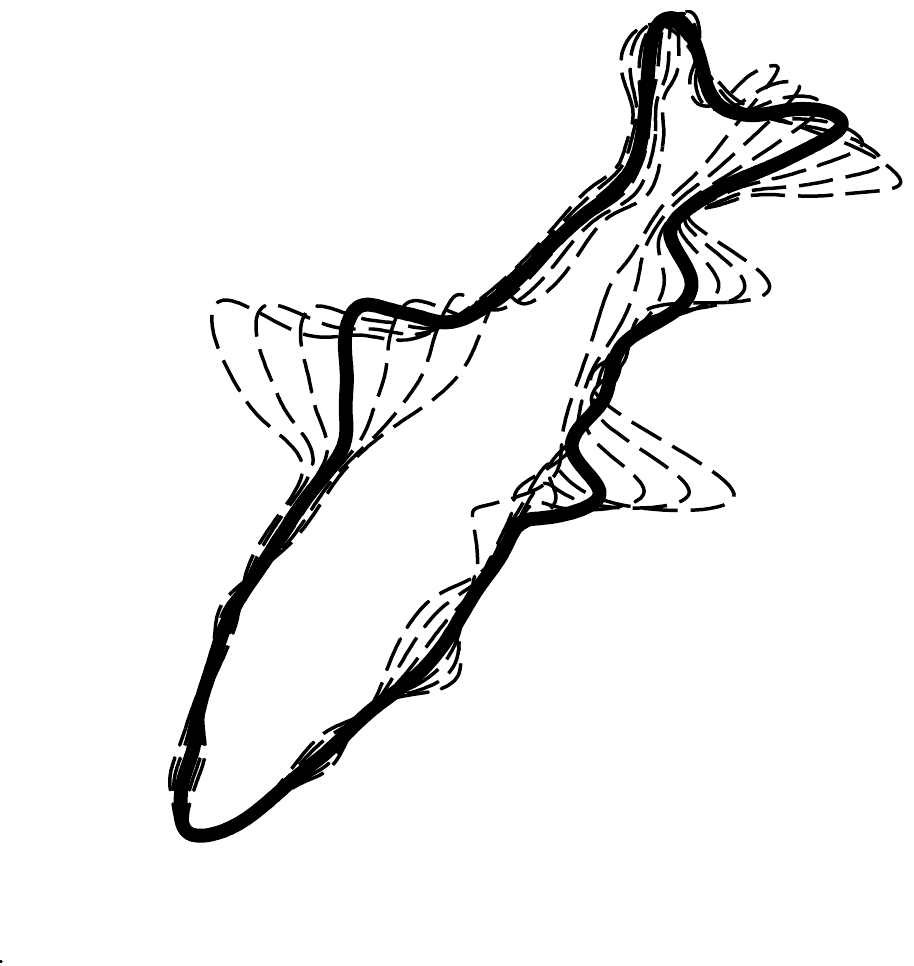}
	\includegraphics[width=.20\textwidth,angle=275]{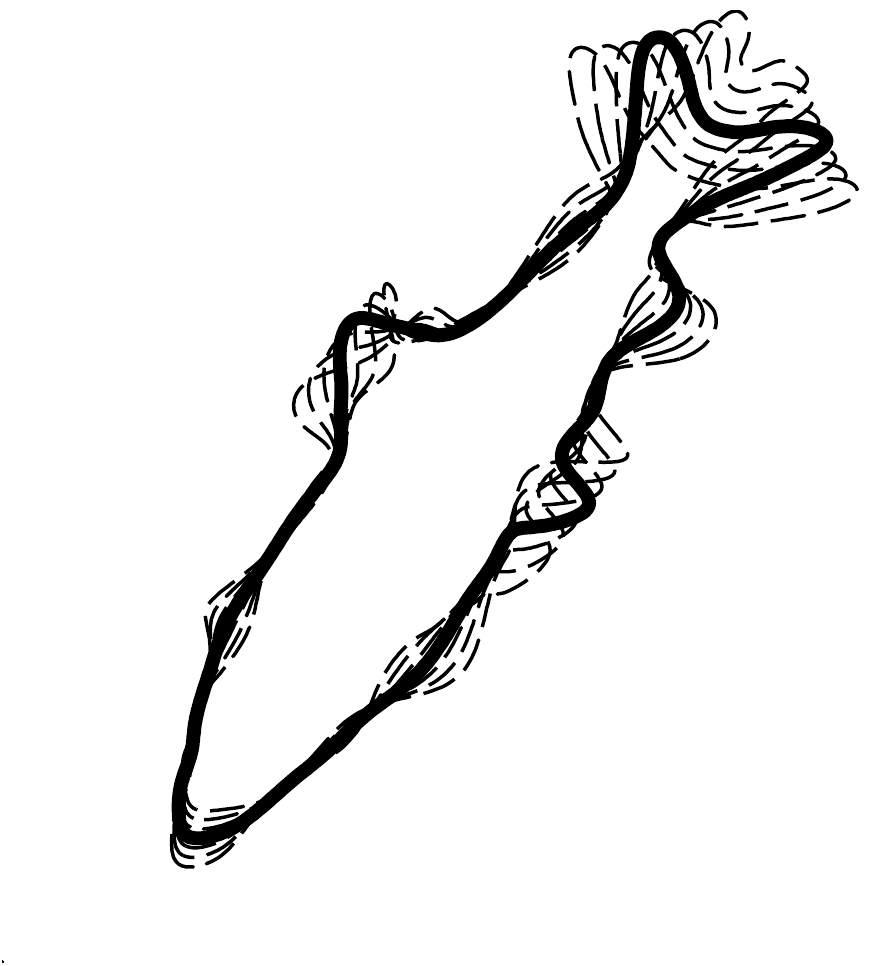}
	\caption{First column: Karcher means (bold) of the groups of fish and humans. Second and third column: geodesics from the mean in the first and second principal direction at the times $-3,-2, \dots, 2, 3$; the bold curve is the mean.}
	\label{fig:KarcherMean}
\end{figure}

Finally, we studied within-group variations using non-linear principal component analysis. To this aim, we first computed the Karcher mean of each group. The corresponding optimization problem \eqref{eq: karcher functional} was solved using a conjugate gradient method, as implemented in the Manopt library \cite{Manopt2014}, on the finite-dimensional spline approximation of the Riemannian manifold of curves. The mean shapes for the groups of fish and humans can be seen in Fig.~\ref{fig:KarcherMean}. Next, we represented each shape in the group by the initial velocity from the mean $\overline{c}$ using the inverse of the Riemannian exponential map. We then performed a principal component analysis with respect to the inner product $G_{\overline{c}}$ on the set of initial velocities. In the group of human figures, the first three eigenvalues capture 67\%, 22\%, and 6\% of within-group variation. In the group of fish, the first three eigenvalues capture only 40\%, 25\%, and 16\% of within-group variation. Geodesics from the mean in the directions of the first two principal directions can be seen in Fig.~\ref{fig:KarcherMean}. In the group of humans the first principal direction encodes bending of the arms and legs, whereas the second direction reflects stretching in the extremities.

\section{Conclusions}
In this article we developed a numerical framework for solving the initial and boundary value problem for geodesics of planar, unparametrized curves under second order Sobolev metrics. We tested our implementation on a dataset of shapes representing various groups of similar physical objects and obtained good experimental results. In future work, we plan to apply our algorithms to datasets of real-world medical data, prove rigorous convergence results for the discretization, and extend the framework to other spaces of mappings like manifold-valued curves and embedded surfaces.

\section*{Acknowledgements}
We want to thank Peter W. Michor and Jens Gravesen for helpful discussions and valuable comments. All authors have been supported by the programme ``Infinite-Dimensional Riemannian Geometry with Applications to Image Matching and Shape Analysis'' held at the Erwin Schr\"odinger Institute. M. Bauer was supported by the European Research Council (ERC), within the project 306445 (Isoperimetric Inequalities and Integral Geometry) and by the FWF-project P24625 (Geometry of Shape spaces).

\end{document}